\documentclass[10pt]{article}

\usepackage{geometry}

\usepackage[utf8]{inputenc} 
\usepackage[T1]{fontenc}    
\usepackage{hyperref}       
\usepackage{url}            
\usepackage{booktabs, amsmath}       
\usepackage{amsfonts}       
\usepackage{nicefrac}       
\usepackage{microtype}      
\usepackage{microtype}
\usepackage{graphicx}
\usepackage{subfigure}
\usepackage{booktabs} 
\usepackage{hyperref}

\usepackage[cmintegrals]{newtxmath}
\usepackage{microtype}
\usepackage{graphicx}
\usepackage{subfigure}
\usepackage{booktabs} 
\usepackage{amssymb,mathrsfs}
\usepackage{bm}
\usepackage{bbm}
\usepackage{subfigure}
\usepackage{graphicx}
\usepackage{calc}
\usepackage{epstopdf}
\usepackage{algorithm,algorithmic}
\usepackage{enumerate,enumitem} 
\newtheorem{thm}{Theorem}[section]
\newtheorem{lemma}{Lemma}[section]

\newtheorem{definition}{Definition}

\usepackage{verbatim}
\usepackage{caption}
\usepackage{multicol,caption}
\usepackage{graphicx}
\usepackage{lipsum}

\usepackage{authblk}
\newcommand{\norm}[1]{\left\lVert#1\right\rVert}

\title{Local Differential Privacy in Decentralized Optimization}
\author[1]{Hanshen Xiao \thanks{hsxiao@mit.edu}}
\author[2]{Yu Ye \thanks{yu9@kth.se}}
\author[1]{Srinivas Devadas \thanks{devadas@csail.mit.edu}}
\affil[1]{CSAIL MIT}
\affil[2]{School of Electrical Engineering and Computer Science, KTH}

\begin{document}



%

	\maketitle


	\begin{abstract}
		Privacy concerns with sensitive data are receiving increasing attention. In this paper, we study local differential privacy (LDP) in interactive decentralized optimization. By constructing random local aggregators, we propose a framework to amplify LDP by a constant. We take Alternating Direction Method of Multipliers (ADMM), and decentralized gradient descent as two concrete examples, where experiments support our theory. In an asymptotic view, we address the following question: {\em Under LDP, is it possible to design a distributed private minimizer for arbitrary closed convex constraints with utility loss not explicitly dependent on dimensionality?} As an affiliated result, we also show that with merely linear secret sharing, information theoretic privacy is achievable for bounded colluding agents. 
	\end{abstract}
	
	\section{Introduction}
	\label{intro}
	Due to the underlying intensive computation and memory requirement in large-scale machine learning, distributed learning has witnessed tremendous development in recent years. In general, there exist two typical scenarios of distributed optimization. The first one assumes a central server to collect and average out local estimates from each agent to update the global model, for example, federated learning \cite{Federated2016},\cite{Federated2017},\cite{cpSGD},\cite{DPML2011}.
	When a "data fusion" center is costly or infeasible, one recourse is a decentralized approach where each agent broadcasts updates to its neighbors and agents collaboratively approach the global optimum \cite{graph1},\cite{graph2},\cite{GD2016},\cite{GD2009}. 
	
	While there is great interest in advances to accelerate the performance of optimization algorithms, privacy preservation is also equally important to many machine learning tasks, especially in the processing of medical records and financial data. Techniques are required to quantify privacy loss during processing, and differential privacy (DP) is one of the best known rigorous theoretical mechanisms that serves this purpose. There is a large body of DP based Empirical Risk Minimization (ERM) work \cite{DPML2011},\cite{ERM2013},\cite{FOCS2014},\cite{lasso}, \cite{Wang2018}, \cite{Wang2017},\cite{ERM2019}, \cite{accuracy2017} and DP-based $k$-means\cite{kmeans}, Bayesian learning \cite{bayesian2017}, identity testing \cite{it} and deep learning \cite{deeplearning2015}, \cite{deeplearning2016} works.
	Given a randomized algorithm, DP offers a provable guarantee against statistical inference that its output is insensitive to a slight change in the input dataset, for example, replacement of a single datapoint. Thus, from outputs observed, it is hard to distinguish the participation of an individual.
	The notion of DP was initially developed with a centralized view.
	Under central DP or distributed learning with trusted central servers, using secure aggregator techniques \cite{distress}, \cite{MPC2015} or subsampling techniques with secrecy of intermediate computation \cite{subsampling},\cite{FOCS2018}, many elegant mechanisms of privacy amplification have been proposed.
	
	However, agents may not trust any other parties to collect their local data. To this end, a stronger notion is local differential privacy (LDP) \cite{localDP}, \cite{local2014} in the distributed scenario, where each agent can run a randomization procedure locally and the privacy of an individual is still guaranteed even with a malicious collector. LDP has been adopted by Apple, Microsoft and Google as one formal definition of privacy \cite{Apple}, \cite{Microsoft}, \cite{Rappor}. 
	Nonetheless, in contrast to the central model, LDP is far less studied despite successful deployment in industry \cite{local2018}. Especially for private optimization, in the central model, the optimization protocol can be viewed as a black box since no information leakage can occur during the execution. Perturbation can be elegantly added only in the objective function, at the beginning, or the output, at the end. Well-known objective/output perturbation methods (e.g. \cite{DPML2011}) and follow-up work \cite{highdimension2012}, \cite{highdimension2014}, \cite{highdimension2016} 
	all capture this idea. For LDP, though one can still apply the above techniques in an non-interactive manner by directly processing the local optimum from each user \cite{SP2017},\cite{Wang2018},\cite{wang2019}, high sample complexity to compensate for accuracy loss may be required. 
	
	In an interactive decentralized optimization, additional privacy loss arises from agents' cooperation. To investigate LDP in general, we need concrete algorithms to support our analysis. In general, there are two types of decentralized optimization. One is (sub)gradient based, such as decentralized (stochastic) gradient descent (GD) methods \cite{GD2009}, \cite{GD2016}, \cite{SGD2017}, and EXTRA \cite{extra}. The second relies on solving a constrained problem with dual variables to minimize some Lagrangian function, such as Alternative Direction Method of Multipliers (ADMM) \cite{convergence}. Though both proceed in an iterative manner, the computation of GD in each step can be less expensive compared to ADMM. Nevertheless, for general convex problems, the convergence rate of decentralized GD is $O(1/\sqrt{K})$ and that of ADMM is $O(1/K)$ \cite{graph1}, where $K$ denotes the number of iterations. Under such a framework, agents enrolled in computing only need to share the states of optimization with neighbors. However, privacy loss also arises from such information exchange, since exposed intermediate results can be easily used to learn the sensitive parameters of the local private functions. Incorporating cryptographic methods, such as (partial) homomorphic encryption \cite{Homoencryt}, \cite{Homoencryt2010}, \cite{IEEE2019} can come with high overhead especially in large-scale optimization. Alternatively, under the lens of DP, the most common approach is to apply perturbed local estimates during the update exchange in decentralized algorithms \cite{ADMM-DP}, \cite{ICML2018}, \cite{DP1}, \cite{SGD2015}, \cite{privateGD2017}, \cite{SGD2018}, \cite{optimal2019}.
	Although heuristic exploration, such as gradually increasing the step penalty \cite{ICML2018}, can improve the utility-privacy tradeoff, existing works lack insights into the fundamentals of algorithmic convergence with noise perturbation. 
	
	\textbf{Problem Statement}: Consider a decentralized optimization problem across $N$ agents in a connected network. The network is modeled by an undirected graph $\mathscr{G}(\mathscr{N}, \mathscr{E})$. Nodes are indexed as $\mathscr{N}=\{1,...,N\}$ and when two nodes $i$ and $j$ are neighbors that can communicate, $(i,j) \in \mathscr{E}$. In general, we assume each node holds a function $f_i(\bm{x}_i)$ that we regard as a loss function determined by samples held locally with the parameter $\bm{x_i}$ to be optimized. Throughout the rest of the paper, we always assume that $f_i(\cdot)$ is a differentiable convex function $\mathcal{C} \to \mathbb{R}$ and $\bm{x}_i \in \mathcal{C} \subset \mathbb{R}^{d}$. $\mathcal{C}$ can be viewed as the constraint, assumed to be a closed convex set. In general, we express the objective function to minimize as 
	\begin{equation}
	\label{problem}
	\min_{\bm{x}_{[1:N]}}~\sum_{i=1}^{N}f_i(\bm{x}_i),~~ s.t. \sum_{i=1}^{N} A_i \bm{x}_{i} = \bm{c}, 
	\end{equation}
	under a linear constraint. In many learning problems, $\bm{x}_{[1:N]}$ stand for one parameter to be collaboratively optimized, where $[1:N]$ is the compact form of $\{1,2,...,N\}$. We term the problem as consensus optimization if the constraint requires that all $\bm{x}_i$ be equal, which can also be rewritten as a linear constraint $\sum_{i=1}^{N} A_i \bm{x}_{i} = \bm{0}$, where $A_i$ includes the graph connectivity \cite{graph1}.
	
	\textbf{Differential Privacy}: For a randomized algorithm $\mathscr{A}$ and a dataset $\mathscr{D}$ as its input, we call $\mathscr{D'}$ adjacent to $\mathscr{D}$ if $\mathscr{D}$ and $\mathscr{D'}$ only differ in one data point.
	Quantitatively, we say $\mathscr{A}$ achieves $\epsilon$-DP if for any $\mathscr{D'}$ and any set $S$ in the domain of $\mathscr{A}(\cdot)$, $ \sup_{\mathscr{D'}, S}  | \log(\Pr[\mathscr{A}(\mathscr{D}) \in S]) - \log( \Pr[\mathscr{A}(\mathscr{D'}) \in S]) | \leq \epsilon$.
	Relaxed DP in the form $(\epsilon, \delta)$ can be found in \cite{relaxedDP}. To embed the notion of DP in the decentralized optimization setting, in this paper, $\mathscr{A}$ corresponds to the optimization algorithm selected, while the functions $f_i$ behave as the inputs and are the privacy concern.
	In the local version, each agent does not trust anyone and, in the worst case, all other parties are colluding against some agent $i$ to learn some sensitive information of $f_i$. With symmetry, we consider the view of agent $i$ and select $\mathscr{D}=\{f_1,...,{f}_i,...,f_N\}$ and $\mathscr{D'}=\{f_1,...,\hat{f}_i,...,f_N\}$, where $\hat{f}_i$ represents some possible candidate of $f_i$.
	Following \cite{localDP}, we define a quantifiable measure, \textbf{$\epsilon$-LDP}, of $\mathscr{A}$ on some output $\bm{S}$ as
	\begin{equation}
	\label{localdp_basic}
	\boxed{\epsilon(\bm{S}) = \sup_{\mathscr{D}, \mathscr{D'}}  \big|\log(\Pr[\mathscr{D} |\bm{S}]) - \log(\Pr[\mathscr{D'} |\bm{S}] ) \big| }  
	\end{equation}
	Notice that we have dropped the dependence of $\epsilon(S)$ on the index $i$ for notational brevity.
	
	\textbf{Methodology Overview}: We first review the basic updating subroutine in decentralized gradient descent. Let $\bm{x}^{k}_i$ denote the local estimate of the global optimum for agent $i$ at round $k$. In consensus optimization, a naively noisy form of decentralized GD \cite{GD2009}, \cite{GD2016} can be described as,
	\begin{equation}
	\label{GD_rule}
	\boxed{ \bm{x}^{k+1}_i := \underline{\sum_{i=1}^{N} w_{ij}\bm{x}^k_j}_{A} -\underline{\eta_{k+1}\nabla f_i(\bm{x}^k_i)}_{B}+\bm{\Delta}^{k+1}_i} 
	\end{equation}
	where $w_{ij} \in [0,1]$ is the weight assigned to $\bm{x}^k_j$ such that $\sum_{j=1}^{N} w_{ij}=1$ and $\eta_{k+1}$ is the step size at round $(k+1)$.
	$\bm{\Delta}^{k+1}_i$ denotes the noise added, which is assumed to be Laplace noise. Therefore, Part A in (\ref{GD_rule}) is a weighted average of updates collected from the last iteration and Part B corresponds to the gradient step taken. 
	In the typical case (\ref{GD_rule}), gradient descent is randomized by resorting to perturbation $\bm{\Delta}^{k+1}_i$. One can apply a stochastic gradient instead, but the subsampling does not amplify privacy \cite{FOCS2014}. However, we show that carefully designed randomization can reshape the output distribution to produce a sharpened privacy loss. In Section 2, it is noted that the optimization protocol may take dozens of steps and divergence among $\bm{x}^{k}_{i}$, $i \in [1:N]$, always exists. Now, for instance, imagine if we use random weights $\widetilde{w}_{ij}$ instead in (\ref{GD_rule}), which are independently selected across each iteration: Part A becomes a random variable within some neighborhood of $\bm{x}^{k}_{i}$, while Part B remains the same.
	In the $(k+1)$th iteration, conditional on $\bm{x}^{k}_{i}$, $i \in [1:N]$, and $f_i$, $\bm{x}^{k+1}_i$ then follows a mixture Laplace distribution with a random mean. Similarly, recalling the updating rule of ADMM with dual method for (\ref{problem}), each update relies on solving an optimization problem:
	\begin{equation}
	\label{leak1}
	\boxed{ \bm{x}^{k+1}_i := \arg \min_{\bm{x}_i} \mathcal{L} (\bm{x}^k_1, ... ,\bm{x}_i,...,\bm{x}^k_N, \bm{\lambda}^k)  + \frac{\rho}{2} \bigg\| A_i \bm{x}_i + \sum_{j \not=i}^{N} A_j \bm{x}^k_j -\bm{c}\bigg\|^2 + \frac{\Gamma}{2} \big\|\bm{x}_i - \bm{x}^k_i\big\|^2+\bm{\Delta}^{k+1}_i} 
	\end{equation}
	where the Lagrangian function is defined as $ \mathcal{L} (\bm{x}_1, ... ,\bm{x}_N, \bm{\lambda})= \sum_{i=1}^{N} f_i(\bm{x}_i) - \bm{\lambda}^T(\sum_{i=1}^{N} A_i \bm{x}_i-\bm{c} ).$
	The Lagrangian multiplier $ \bm{\lambda}^{k+1}$ is updated through $	\bm{\lambda}^{k+1} := \bm{\lambda}^{k} - \zeta (\sum_{i=1}^N A_i\bm{x}^{k+1}_i -\bm{c}).$ If we allow each agent to independently select random penalties $\rho$ and $\Gamma$ across iterations, clearly, a similar mixture Laplace with random mean, $\bm{x}^{k+1}_i$ is produced. While a random mean amplifies the uncertainty, trivially incorporating such an idea in existing algorithms may easily lead to unclear convergence guarantees or cumbersome privacy analysis.
	
	\textbf{Contribution and organization}: In this paper, we study the LDP of interactive decentralized optimization from both asymptotic and non-asymptotic perspectives. 
	\begin{enumerate}[label=(\roman*)]
		\item Taking GD and ADMM as two examples, we propose a framework of decentralized optimization with varying parameters. A unified local privacy analysis is presented in Theorem 2.2 with assumptions only on local sensitivity and we quantify the privacy amplification constant in Theorem 2.3. Rigorous analysis on the convergence rate and the utility-privacy tradeoff for proposed algorithms are included in Theorems 3.1-3.3 and Theorem 3.4, for ADMM and GD, respectively. Appealingly, in contrast to noise perturbation, which always brings a compromise in utility, we show with carefully selected parameters, the proposed schemes enjoy the same convergence rate as before. Experiments support the theory.
		
		\item From an asymptotic viewpoint, we follow a rich line of work to study LDP in high dimensional optimization. In the central model, it is well understood that with certain sparsity of $\mathcal{C}$, utility loss can scale logarithmically in dimensionality $d$, for example in sparse linear regression (LASSO) \cite{lasso}. In this paper, we explore such dimensional dependence from the objective function $f_i$'s side. Motivated by \cite{highdimension2014} we show even under LDP, either for generalized linear functions (Theorem 4.1), or when $\nabla f_i$ is bounded in $l_1$ norm (Theorem 4.2), the utility loss is not explicitly dependent on (ignoring logarithmic terms)  $d$ for arbitrary closed convex constraints $\mathcal{C}$. 
	\end{enumerate}

	\section{Algorithm Description and Privacy Analysis}
	We first present the main protocol of modified private ADMM and GD as Algorithm 1 and 2, respectively. Here, we omit the projection step if $ \mathcal{C} \subsetneqq \mathbb{R}^d$ since we are interested in arbitrary constraints. For simplicity, in this section we assume $\mathscr{G}$ is fully connected and only consider the consensus problem temporarily, while we provide a convergence proof for general cases in the next section. In contrast to previous private ADMM protocols \cite{convergence}, \cite{ADMM-DP}, \cite{ICML2018}, we consider applying a first-order approximation for each $f_i$, similar to \cite{f2011}, \cite{f2015}, as
	\begin{equation}
	\label{firstorder}
	f_i(\bm{x}_i) \approx f_i(\bm{x}^k_i) + \nabla f_i(\bm{x}^k_i) (\bm{x}_i-\bm{x}^k_i),
	\end{equation}
	and we derive Algorithm 1 accordingly. Thus, Algorithm 1 has reduced complexity, same as that of GD, which avoids solving (\ref{leak1}). In this paper, we use $\norm{\cdot}_q$ to denote the $l_q$ norm and $\norm{\cdot}$ denotes the standard $l_2$ norm for brevity. It is clear that Algorithm 1 and 2 share a very similar structure except for the dual variable $\bm{\lambda}^k_i$ in ADMM. Intuitively, Term (A) in either (\ref{perturbation-updating}) or (\ref{GD_new}) behaves as a random aggregator to merge the updates from the previous iteration and Part $B$ corresponds to the effect from the function $f_i$ on updating $\bm{x}^{k+1}_i$. Here, we only give two examples to select the random weights where Term (A) is a uniformly distributed in some interval. Many variants can be easily derived, and more details can be found in Appendix E. Clearly, in Algorithms 1 and 2, only the exchange of $\bm{x}^{k}_i$ incurs privacy leakage. Hence, if secure aggregation exists in a decentralized manner, information-theoretic privacy is achievable. 
	\begin{thm}
		\label{info-sharing-thm}
		If each node can have secure communication with $(N-1)$ nodes, ADMM/GD achieves information-theoretic privacy with no utility compromise via linear secret sharing, when there are at most $(N-2)$ colluding nodes, $N \geq 3$. 
	\end{thm}
	The proof can be found in Appendix A. Though there is no compromise in accuracy, the limitation of the above theorem is also clear as it only works for bounded colluding agents. Also, there is no guarantee associated with how much information can be learned from the convergence. Recalling Algorithm 1 and 2, let $\bm{X} = \bm{x}^{[0:K]}_{[1:N]}$ denote the outputs across $(K+1)$ iterations. The definition of LDP for agent $i$ is formally given as follows. 
	\begin{algorithm}[tb]
		\caption{Modified Private ADMM with First-order Approximation}
		\label{new-admm}
		\begin{algorithmic}
			\STATE {\bfseries Input:} Local functions $f_{[1:N]}$, step penalty $\zeta$. 
			\STATE Initialize $\bm{x}^0_{[1:N]}$ randomly, $\bm{\lambda}^0_{[1:N]} =\bm{0}$. Each agent selects a private constant $D_i$. 
			\FOR{$k=0,1,2, ... K-1$}
			\STATE {\bfseries Agents $i=1$ to $N$} do in parallel:\\
			\STATE Randomly pick two positive diagonal matrices $\bar{\bm{\rho}}^{k+1}_i$ and $\bm{\Gamma}^{k+1}_i$ such that $(N-1)\bar{\bm{\rho}}^{k+1}_i + \bm{\Gamma}^{k+1}_i = D_i \cdot \bm{I}_{p}$ and update $\bm{x}^{k+1}_i$ \textit{in parallel}:
			\begin{equation}
			\label{perturbation-updating}
			\bm{x}^{k+1}_i :=\underline{\frac{\bm{\Gamma}_i^{k+1}}{D_i}\bm{x}_i^k + \frac{(N-1)\bar{\bm{\rho}}_i^{k+1}}{D_i} \frac{\sum_{ j \not = i}\bm{x}_j^k}{N-1}}_A - \underline{D^{-1}_i\nabla f_i(\bm{x}^k_i)}_B + D^{-1}_i\bm{\lambda}^{k}_i +\bm{\Delta}^{k+1}_i .
			\end{equation}
			\STATE Exchange $\bm{x}^{k+1}_i$ and then update $\bm{\lambda}^{k+1}_i := \bm{\lambda}^{k}_i - \zeta \sum_{j \not = i} (\bm{x}^{k+1}_j -\bm{x}^{k+1}_i)$.
			\ENDFOR
		\end{algorithmic}
	\end{algorithm}
	
		\begin{algorithm}[tb]
		\caption{Modified Private Decentralized Gradient Descent}
		\label{new-GD}
		\begin{algorithmic}
			\STATE {\bfseries Input:} Local functions $f_i$ and a diminishing sequence $\{\eta_{k}\}$
			\STATE Initialize $\bm{x}^0_{[1:N]}$.  
			\FOR{$k=0,1,2, ...,K-1 $}
			\STATE {\bfseries Agents $i=1$ to $N$} do in parallel :\\
			\STATE Randomly pick a positive diagonal matrix $\bm{w}^{k+1}_{i}$ of which the non-zero elements are within $(0,1)$. Then, update $\bm{x}^{k+1}_i$:
			\begin{equation}
			\label{GD_new}
			\bm{x}^{k+1}_i := \underline{ \frac{2\bm{w}^{k+1}_{i}}{N} \bm{x}^k_{i_k}+ \frac{2(\bm{I}_{d}-\bm{w}^{k+1}_{i})}{N} \bm{x}^k_{\hat{i}_k} + \frac{\sum_{j \not = i_k, \hat{i}_k} \bm{x}^k_j }{N}}_{A} -\underline{\eta_{k+1}\nabla f_i\big(\frac{\sum_{i=1}^{N} \bm{x}^k_i}{N} \big)}_{B}+\bm{\Delta}^{k+1}_i,
			\end{equation}
			where $(i_k, \hat{i}_k) =\arg \max_{(j_1,j_2)} \norm{\bm{x}^{k}_{j_1}- \bm{x}^{k}_{j_2}}_1$. Exchange $\bm{x}^{k+1}_i $ .
			\ENDFOR
		\end{algorithmic}
	\end{algorithm}
	
	

	
	
	
	\begin{definition} [$(\epsilon(\bm{X}), \delta(\bm{X}))$-LDP] \cite{NISP2018}, \cite{optimal2019}: A decentralized optimization $\mathscr{M}$ achieves $(\epsilon(\bm{X}),\delta(\bm{X}))$-LDP on an output $\bm{X}$, if for any two adjacent inputs $\mathscr{D}=\{f_1, ... ,f_i, ... ,f_N\}$ and $\mathscr{D'}=\{f_1, ... ,\hat{f}_i, ... ,f_N\}$, which differ only in $f_i$ and $\hat{f}_i \in \mathscr{F}_i$, 
		\begin{equation}
		\sup_{\mathscr{D},\mathscr{D'}}  {\Pr[ \mathscr{D}  | \bm{X} ]} \leq \epsilon(\bm{X}) { \Pr[  \mathscr{D'} | \bm{X} ]} + \delta(\bm{X}), 
		\end{equation}
		$\mathscr{F}_i$ is some set of functions.
	\end{definition}
	
	When $\delta(\bm{X})=0$ predefined, we term it pure LDP, otherwise we call it relaxed or approximate LDP, where one may apply a stronger composition theorem for accumulated privacy loss \cite{relaxedDP}. We first focus on the pure LDP case. To quantify the privacy loss, we need to introduce a concrete assumption, {\em{sensitivity}}, to further define $\mathscr{F}_i$. We say $\mathscr{F}_i$ with $\mathscr{B}_q$ local sensitivity if for any $\hat{f}_i, f_i \in \mathscr{F}_i$, $\mathscr{B}_q \geq \sup_{\bm{x} \in \mathcal{C}} \norm{f_i(\bm{x})-\hat{f}_i(\bm{x})}_{q}$ in $l_q$ norm. For example, $\mathscr{B}_{\infty}$ captures the case of an ERM whose data set is bounded in $l_{\infty}$. 
	In the following, we restrict our focus to $\mathscr{F}_i$ with $\mathscr{B}_{\infty}$ local sensitivity. The following lemma, whose proof is in Appendix B, provides a semi closed-form of $\epsilon(\bm{X})$. 
	\begin{lemma} 
		Algorithm 1 and 2 satisfy $\epsilon(\bm{X})$-LDP, where 
		\begin{equation}
		\epsilon(\bm{X}) = \sup_{f_i, \hat{f}_i \in \mathscr{F}_i}\left | \sum_{k=0}^{K-1} \sum_{l=1}^d \log \frac{P\big( \bm{x}^{k+1}_{i}[l] \big| f_i, \bm{x}^{[0:k]}_{[1:N]}\big)}{ P\big( \bm{x}^{k+1}_{i}[l]  \big| \hat{f}_i, \bm{x}^{[0:k]}_{[1:N]}\big)}  \right |, 
		\end{equation}
		where $\bm{x}^{k}_{i}[l]$ denotes the $l^{th}$ coordinate of $\bm{x}^{k}_{i} \in \mathbb{R}^d$. 
	\end{lemma}
	Following the Laplace mechanism \cite{laplace},\cite{laplace2}, we assume each coordinate of $\bm{\Delta}_i^{k}$ i.i.d. following a Laplace distribution Lap$(0,\beta_k)$ with probability density $P(y) = \frac{\beta_k}{2}e^{-\beta_k|y|}$. Thus, either in (\ref{GD_new}) or (\ref{perturbation-updating}), $P( \bm{x}^{k+1}_{i}[l] | f_i, \bm{x}^{[0:k]}_{[1:N]})$ follows the same Laplace distribution as that of $\bm{\Delta}_i^{k+1}[l]$ except the mean of $\bm{x}^{k+1}_{i}[l]$ is a random variable uniformly distributed in some interval, denoted by $\tau^{k+1}_i[l]$ (see Term (A) in (\ref{perturbation-updating}) and (\ref{GD_new}), respectively). To bound the accumulated privacy loss, we consider the composition of the loss on each coordinate $l$ and iteration $k$, denoted by ${\epsilon}^{k+1}_l(\bm{X}) = \sup_{f_i, \hat{f}_i} \left | \log[P( \bm{x}^{k+1}_{i}[l] | f_i, \bm{x}^{[0:k]}_{[1:N]} / {P( \bm{x}^{k+1}_{i}[l] | \hat{f}_i, \bm{x}^{[0:k]}_{[1:N]})}]  \right |.$ The following theorem provides an upper bound of ${\epsilon}^{k+1}_l(\bm{X})$.

	\begin{thm} Algorithms 1 and 2 achieve $\epsilon(\bm{X}) \leq \sum_{k=0}^{K-1} \sum_{l=1}^{d} {\epsilon}^{k+1}_l(\bm{X})$-LDP, \footnote{One can obtain a stronger composition form in relaxed $(\epsilon,\delta)$-LDP \cite{relaxedDP} that Algorithms 1 and 2 achieve $(\sum_{k=0}^{K} \sum_{l=1}^{d} \frac{(e^{\epsilon^{k+1}_l(\bm{X})}-1)\epsilon^{k+1}_l(\bm{X})}{e^{\epsilon^{k+1}_l(\bm{X})}+1} + \sqrt{-2\sum_{k=0}^{K} \sum_{l=1}^{d} (\epsilon^{k+1}_l(\bm{X}))^2 \log(\delta)}, \delta)$-LDP, for some $\delta \in (0,1).$} where 
		\begin{equation}
		\label{worst-upper-bound}
		{\epsilon}^{k}_l(\bm{X})  \leq \max_{|t| \leq \alpha_k\mathscr{B}_{\infty} } \bigg| \log\bigg[\int_{ \tau^{k}_i[l] } e^{-{\beta_k}{\left| \bm{x}^{k}_{i}[l]-x\right| }} dx\bigg] -  \log\bigg[{\int_{ \tau^{k}_i[l]+t } e^{-{\beta_k}{\left| \bm{x}^{k}_{i}[l]-x\right| }} dx }\bigg]  \bigg|,
		\end{equation}
		where $\alpha_k=D^{-1}_i$ in Algorithm 1 while $\alpha_k = \eta_{k}$ in Algorithm 2. $\tau^{k}_i[l]+t$ implies uniformly moving the interval with $t$. Moreover, for arbitrary $\bm{X}$, the right hand of (\ref{worst-upper-bound}) is never bigger than $\alpha_k\beta_k\mathscr{B}_{\infty}$. Specifically, when $\bm{x}^{k}_{i}[l]$ belongs to $\tau^{k}_i[l]$, it is strictly smaller than $\alpha_k\beta_k\mathscr{B}_{\infty}$. 
	\end{thm}
	As a straightforward corollary of the above theorem, if we fix parameters in Algorithms 1 and 2, then the local privacy loss is $\sum_{k=0}^{K-1} d\alpha_{k+1}\beta_{k+1}\mathscr{B}_{\infty}$, which matches prior results \cite{ICML2018}, \cite{ADMM-DP}, \cite{optimal2019}. We include the proof of the above theorem in Appendix C. To conclude, for each $\epsilon^k_l({\bm{X}})$, randomized weights renders a constant privacy reduction expressed as a conditional expectation, $\gamma^k_l = \mathbb{E}_{\bm{x}^k_i}[\epsilon^k_l({\bm{X}})/(\alpha_k \beta_k\mathscr{B}_{\infty})| \bm{x}^{[0:k-1]}_{[1:N]}]$. From (\ref{worst-upper-bound}), a longer length of the interval $\tau^k_i[l]$, denoted by $\omega$, renders a more concentrated mixture distribution compared to pure Laplace. In the following theorem, whose proof is in Appendix D, we quantify $\gamma^k_l$ when $\omega>\alpha_k\mathscr{B}_{\infty}$.
	\begin{thm}
		\label{advantage}
		When $\omega > \alpha_k\mathscr{B}_{\infty}$,
		\begin{equation}
		\gamma^k_l \leq \frac{1}{{\alpha^k\beta_k\mathscr{B}_{\infty}}}\log \bigg\{ e^{\alpha_k\beta_k\mathscr{B}_{\infty}} \bigg[1-2\int_{0}^{\frac{\omega-\alpha_k\mathscr{B}_{\infty}}{2}}\int_{0}^{\omega} \Phi(x,y) dy dx\bigg] + 2 \int_{0}^{\frac{\omega-\alpha_k\mathscr{B}_{\infty}}{2}} \int_{-\alpha_k\mathscr{B}_{\infty}}^{\omega-\alpha_k\mathscr{B}_{\infty}}\Phi(x,y) dy dx \bigg\},
		\end{equation}
		where $\Phi(x,y)=\frac{\beta_k}{2\omega}e^{-{\beta_k}{| x -y|}}.$
	\end{thm}
	In Fig. 1, we show the relationship between $\gamma$, $\omega$ and $\beta$, where temporarily the dependence on $l$ and $k$ is dropped for brevity and $\alpha\mathscr{B}_{\infty}$ is fixed to $0.001$. With $\omega$ ranging from $0.1$ to $1$ and $\beta$ from $2$ to $10$, clearly larger $\omega$ and $\beta$, corresponding to a longer interval length and noise of smaller variance, lead to better privacy amplification. 
	\begin{figure}[h!]
		\centering
		\includegraphics[width=1 \linewidth]{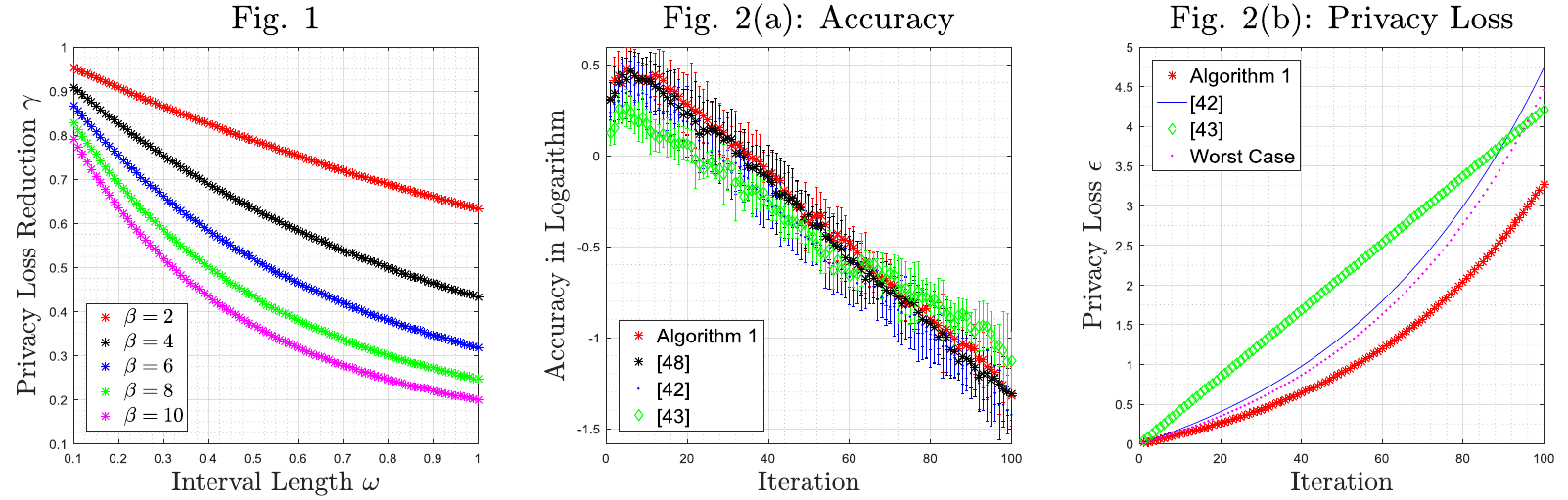}
	\end{figure}
	As a comparison, we also test the proposed schemes and state-of-the-art approaches on a regularized logistic regression of the $Adult$ dataset, from the UCI machine learning repository \cite{ICML2018}.
	The performance of existing private ADMM \cite{ICML2018} (with increasing penalty), \cite{ADMM-DP}, \cite{optimal2019} and Algorithm 1 is illustrated in Fig. 2, where the communication graph $\mathscr{G}$ is randomly generated with $N=10$ and $|\mathscr{E}|=20$. Fig. 2 (a) shows the accuracy in a logarithmic scale, where even with the first-order approximation, Algorithm 1 has almost the same performance as prior works. Associated privacy loss of Algorithm 1 is presented in Fig. 2 (b). The worst case in Fig. 2 (b) refers to $\sum_{k=0}^{K-1}d\alpha_{k+1}\beta_{k+1}\mathscr{B}_{\infty}$. On average, we achieve $30\%$ privacy loss reduction. We omit the privacy loss of \cite{optimal2019} in Fig. 2 (b) since it is too loose in this example. We note that the plots of \cite{ICML2018}, \cite{ADMM-DP} in Fig. 2(b) (and \cite{optimal2019}) either require global sensitivity or smoothness of gradients, while we only assume local sensitivity. Full description and results of experiments are included in Appendix E. In Appendix E (Fig. D), we also show that the parameter randomization in Algorithms 1 and 2 does not bring accuracy compromise, which is guaranteed by the convergence theorems presented in next section. 
	
	



	\section{Convergence and Utility Tradeoff Analysis}
	Starting from this section, we provide convergence analysis of the two proposed decentralized algorithms and study the asymptotic upper bound of utility loss. We will drop the dependence of privacy parameters $(\epsilon,\delta)$ on $\bm{X}$ since we are concerned with the worst case. We first review some commonly used concepts in convex optimization. A function $f(\bm{x}): \mathcal{C} \to \mathbb{R}$ is $L$-Lipschitz continuous if for any $\bm{x}, \bm{y} \in \mathcal{C} $, $|f(\bm{x}) - f(\bm{y})| \leq L \norm{\bm{x}-\bm{y}}$. $f(\bm{x})$ is $m$-strongly convex if $m \norm{\bm{x}-\bm{y}}^2 \leq (\bm{x}-\bm{y})^T(\nabla f(\bm{x})- \nabla f(\bm{y}))$. Further, we call $f(\bm{x})$ $M$-smooth if $\nabla f(\bm{x})$ is $M$-Lipschitz continuous, i.e., $\norm{\nabla f(\bm{x}) - \nabla f(\bm{y})} \leq M \norm{\bm{x}-\bm{y}}.$
	
	We provide convergence proofs of proposed algorithms with different assumptions. When $f_i(\bm{x}_i)$ is $m_i$-strongly convex and $\sqrt{M_i}$-smooth, we describe our construction of Algorithm 1 on ADMM in three steps to solve (\ref{problem}) in general. First, we show the admissible range of random selection of parameters which still preserves a deterministic linear convergence rate (Theorem 3.1). Second, to further reduce the computational overhead, rather than solving equation (\ref{leak1}), we propose a first-order based approximation (Theorem 3.2). In the modified ADMM, the computation in each iteration is simplified to a closed form. Third, we present the hybrid ADMM version with noise perturbation and a utility-privacy tradeoff (Theorem 3.3). 
	To solve (\ref{problem}), with random penalties across iterations, the updating procedure of node $i$ at the $(k+1)$th iteration becomes 
	\begin{equation}
	\label{time-varying-x}
	\bm{x}^{k+1}_i := \arg \min_{\bm{x}_i} f_i(\bm{x}_i) - \bm{\lambda}^{kT}\bigg(A_{i}\bm{x}_i+\sum_{j \not = i} A_j\bm{x}^{k}_j-\bm{c}\bigg)+ \frac{1}{2}  \bigg\|A_{i}\bm{x}_i+\sum_{j \not = i} A_j \bm{x}^{k}_j-\bm{c}\bigg\| ^2_{\bm{\rho}^{k+1}_{i}} + \frac{1}{2} \big\|\bm{x}_i-\bm{x}^{k}_i\big\|^2_{\bm{\Gamma}^{k+1}_i},
	\end{equation}
	and the Lagrangian multiplier is updated accordingly as $\bm{\lambda}^{k+1} := \bm{\lambda}^{k} - \bm{\gamma}^{k+1}_i\bm{\rho}^{k+1}_i(\sum_{i=1}^{N} A_i \bm{x}^{k+1}_i-\bm{c}).$ Here, $\bm{\lambda}^{kT}=(\bm{\lambda}^k)^T$ and $\bm{\gamma}^{k+1}_i\bm{\rho}^{k+1}_i = \zeta \cdot \bm{I}$ is a global constant set up at the beginning and $\norm{\bm{z}}^2_{G} = \bm{z}^TG\bm{z}$. Let $\bm{u}^{k+1} = [\bm{x}^{k+1}_{[1:N]},\bm{\lambda}^{k+1}]$ and $\bm{u}^{*}$ stand for the optimum to (\ref{problem}). 
	\begin{thm}
		\label{varying-rho}
		The proposed ADMM converges linearly to $\bm{u}^{*}$ with penalty $D_i \cdot \bm{I}= A^T_i\bm{\rho}^{k+1}_iA_i + \bm{\Gamma}^{k+1}_i$, where $D_i$ is a constant, if 
		\begin{equation*}
		\begin{aligned}
		\alpha <\frac{2m_i}{N\rho_{i,\max}^{2(k+1)}\sigma_{i,\max}^2+ \breve{\rho}_{i,\max}^{2(k+1)}\sigma^2_{i,\max}},~~  \rho^0 >\frac{N}{2\alpha}, D_i >\max\bigg\{ \rho_{i,\max}^{k+1}\sigma_{i,\max}^2,\frac{N\sigma^2_{i,\max}}{\alpha}  \bigg\},~~\zeta < 2\rho^0 - \frac{N}{\alpha},
		\end{aligned}
		\end{equation*}
		for some positive $\alpha$, $\zeta$ and $\rho^0$. Here $\breve{\rho}_{i,\max}^{ k+1}$ is the diagonal element of matrix $\bm{\rho}_i^{k+1}- \rho ^0\cdot\bm{I}$ with the maximal absolute value and $\sigma_{i,\max}$ is the largest singular value of $A_i$. More specifically, $$\norm{\bm{u}^{k}-\bm{u}^{*}}^2_G \geq (1+p) \norm{\bm{u}^{k+1}-\bm{u}^{*}}^2_G,$$ for some $p>0$, where $G = diag( D_1 \cdot \bm{I}, ... ,D_N \cdot  \bm{I}, \zeta \cdot \bm{I})$. The selection of $p$ is specified in (\ref{delta-selection}) in Appendix F.
	\end{thm}
	From Theorem \ref{varying-rho}, it is noted that both $\bm{\rho}^{k}_{[1:N]}$ and $\bm{\Gamma}^{k}_{[1:N]}$ are not necessarily constant. When ${D_i}$ is sufficiently large and $\zeta$ is sufficiently small, which indicates that the step sizes of both primal variable $\bm{x}_i$ and dual variable $\bm{\lambda}$ are small enough, $\bm{\rho}^{k}_{[1:N]}$ can be independently and randomly selected from an interval centered at some point $\rho^0 \cdot \bm{I}$ and $\bm{\Gamma}^{k+1}_i = D_i \cdot \bm{I} - A^T_i\bm{\rho}^{k+1}_iA_i$. However, with (\ref{time-varying-x}), we may encounter considerable computation overhead in each iteration when no closed-form optima exists. Substituting (\ref{firstorder}) into (\ref{time-varying-x}), $\bm{x}^{k+1}_i$ can then be approximated in a closed-form of $\bm{x}^k_{[1:N]}$ and $\bm{\lambda}^k$,
	\begin{equation}
	\label{first-order-updating}
	\bm{x}^{k+1}_i := \bm{D}^{-1}_i \big[ A^T_i\big( \bm{ \lambda}^{k}-  \bm{\rho}^{k+1}_i \big( \sum_{j \not = i} A_j\bm{x}^{k}_j- \bm{c}\big) \big) + \bm{\Gamma}^{k+1}_i\bm{x}^{k}_i  -  \nabla f_i(\bm{x}^k_i)\big].
	\end{equation}
	To quantify the loss from the approximation, we provide the following theorem.
	\begin{thm}
		\label{first-order-converge}
		First-order approximation based ADMM, with modified updating procedure (\ref{first-order-updating}) still enjoys the linear convergence rate with proper penalty selection specified in (\ref{rho-fa-selection}) in Appendix G.
	\end{thm}
	Now, we further consider the perturbation version of (\ref{first-order-updating}) with noise: the only difference is that an independent noise $\bm{\Delta}^{k}_{i}$ is added at the end of the updating procedure as shown in (\ref{perturbation-updating}). 
	\begin{thm}
		\label{noise-thm}
		With the same assumptions as Theorem \ref{first-order-converge}, if the updating procedure further perturbs with an independent noise $\bm{\Delta}^{k+1}_{i}$, defined in (\ref{perturbation-updating}) in Algorithm 1, there exists a constant $a \in(0,1)$ and residual $R^{k}$ such that 
		\begin{equation}
		\label{deviation-bound}
		\norm{\bm{u}^{k} - \bm{u}^{*}}^2_G \leq a^k \norm{\bm{u}^{0} - \bm{u}^{*}}^2_G+ R^{k},
		\end{equation}
		where the expression of $a$ and $R^{k}$ can be found in Appendix H. \footnote{In particular, when $\lim_{k \to \infty} \bm{\Delta}^{k}_i \to \bm{0}$ for each $i$, i.e., a diminishing noise is utilized, $\lim_{k \to \infty} R^{k} \to 0$.} With a total $\epsilon$ LDP budget, $\mathbb{E} \norm{\bm{u}^{k} - \bm{u}^{*}}^2 = \widetilde{O} (\frac{d^3N\mathscr{B}^2_{\infty}}{\epsilon^2})$ and under relaxed $(\epsilon,\delta)$-LDP, this bound is sharpened to $\widetilde{O}(\frac{d^2N\mathscr{B}^2_{\infty}}{\epsilon^2})$, where we ignore other constants with respect to $A_i$ and $f_i$. Here $\widetilde{O}$ is the big-O that ignores logarithmic factors.
	\end{thm}

	The proof of Theorem 3.3 can be found in Appendix H. In the general convex case, we take GD as an example. In comparison to ADMM, GD only captures consensus optimization. When $\norm{\nabla f_i}$ is bounded, the following theorem shows the privacy-utility tradeoff of Algorithm 2, with proof in Appendix I. Here we assume $\bm{x}^*$ is the optimum to (1) in the consensus case. 
	\begin{thm}
		\label{GD-convergence}
		Assume that $f_i(\bm{x})$ is convex and $\norm{\nabla f_i(\bm{x})}^2$ is bounded by $G^2$ for each $i$. Moreover, let $V^2 = \max_{i,k} \mathbb{E} (\bm{\Delta}^{k}_i / \eta_{k})^2$. When we select the step size $\eta_k = \frac{1}{c\sqrt{k}}$ for some constant $c$, then 
		$$ \sum_{i=1}^N f_i(\tilde{\bm{x}}^{K-1}) -\sum_{i=1}^N f_i(\bm{x}^*) \leq  \frac{c\sqrt{K} \sum_{i=1}^{N}\norm{\bm{x}^0_i-\bm{x}^*}^2 +N c^{-1}(\log{K}+2)\sqrt{K+1}(G^2 +V^2)}{K},$$
		where $\tilde{\bm{x}}^{K-1} =\frac{1}{K}\sum_{k=0}^{K-1}\bar{\bm{x}}^k $ and $\bar{\bm{x}}^k=\frac{1}{N}\sum_{i=1}^{N}\bm{x}^{k}_i .$ In $\epsilon$-LDP, Algorithm 2 has utility loss $\tilde{O} \bigg(\frac{\sqrt{\sum_{i=1}^{N}\norm{\bm{x}^{0}_i-\bm{x}^{*}}^2})\sqrt{N}(G+\frac{d^{3/2}K\mathscr{B}_{\infty}}{\epsilon})}{\sqrt{K}}\bigg)$. With $(\epsilon, \delta)$ relaxation, this bound can be sharpened to $\tilde{O} \bigg(\frac{\sqrt{\sum_{i=1}^{N}\norm{\bm{x}^{0}_i-\bm{x}^{*}}^2}\sqrt{N}d\mathscr{B}_{\infty}}{\epsilon}\bigg)$ when $K$ is sufficiently large. 
	\end{thm}
	
	\section{Utility Loss under LDP in High-dimensional Optimization}
	Over the last decade, a rich line of work has been devoted to studying and even breaking the curse of dimensionality in DP. On the negative side, Smith et al. in \cite{FOCS2014} have shown that even for linear regression with samples bounded in $l_2$ norm, utility loss has a polynomial dependence on $d$ generally speaking. However, the remarkable work by Talwar et al. \cite{talwar2014} and follow-up works \cite{highdimension2016},\cite{peking2017},\cite{wang2019} show that such dependence can be captured by the Gaussian width of the constraint $\mathcal{C}$ being optimized over. In particular, the Gaussian width over an $l_1$ ball scales logarithmically in $d$ and in \cite{lasso}, a nearly-optimal private Frank-Wolfe algorithm is established to solve LASSO. In addition, some data-dependent assumption or privacy concerns can also help \cite{highdimension2012}, \cite{thakurta2013}, \cite{wang2018high}, \cite{wang2019high}. Such evidence provides us confidence that under certain sparsity, {\em even LDP is possible in high-dimensional cases}. Although the required loss with respect to the constraint $\mathcal{C}$, i.e., the domain of $\bm{x}_i$, is fairly well understood, we are interested in whether the structure of $f_i(\cdot)$ can also help. Under a central model, positive results have been shown for generalized linear functions based on objective/output perturbation \cite{highdimension2014}. In such a setup, for simplicity, we assume $
	f_i(\bm{x}_i) = \frac{1}{b_i} \sum_{j=1}^{b_i} \phi(\langle \bm{z}^i_j, \bm{x}_i \rangle, l^i_j)$, $\bm{x}_i \in \mathcal{C}$, where $(\bm{z}^i_j, l^i_j) \in \mathbb{R}^d \times \mathbb{R}$, $j \in [1 : b_i]$, are data points held by agent $i$ and $\phi(\cdot, \cdot)$ is a convex function $\mathbb{R} \times \mathbb{R} \to \mathbb{R}$. Here, we still focus on the consensus optimization, i.e., $\bm{x}_i = \bm{x}_j$, and thus this formulation captures a supervised learning of a linear predictor. For those data points, if $\bm{z}^i_j$ in $l_2$ norm is bounded by $C_2$ for $i \in[1:N]$ and $j \in [1:b_i]$, then the sensitivity can also be reasonably bounded in $l_2$ norm by $\mathscr{B}_2$. Under such an assumption, we show under LDP, dimensionality independent loss can be derived for generalized linear functions by applying Algorithm 2 on $f_i(\cdot)$. The main idea we use here is to control the error from the Lipschitz assumption of $f_i(\cdot)$ rather than analyzing the divergence $\norm{\bm{x}^k_i - \bm{x}^*}$ directly as before.The proof of the following theorem is in Appendix J.
	\begin{thm}
		\label{glf_thm}
		For $f_i(\bm{x}_i)$ defined in the generalized linear form above, we assume $\phi(\cdot, \cdot)$ is an $L$-Lipschitz continuous function, which is $m$-strongly convex and $M$-smooth. For Algorithm 2 to achieve $(\epsilon, \delta)$-LDP budget, the utility loss $\sum_{i=1}^N f_i(\bar{\bm{x}}^K) - \sum_{i=1}^N f_i(\bm{x}^*)$ can be independent of $d$ and upper bounded by $\tilde{O}\big(\frac{\sqrt{M}NC_2L\mathscr{B}_2}{m^{3/2}\epsilon}\big)$.
	\end{thm}
	
	In general, beyond the generalized linear assumption, as an analog to the case where the constraint $\mathcal{C}$ is $l_1$ bounded, it is natural to consider whether the utility loss can be also expressed in a form $\norm{\nabla f_i}_1$ and with a logarithmic dependence on $d$. In the following, we answer this question affirmatively. We prove the following result in Appendix K.
	\begin{thm}
		\label{l_1}
		If $f_i(\bm{x}_i)$ is $m$-strongly convex and $M$ smooth, of which $\norm{\nabla f_i(\bm{x}_i)}_1 \leq C_1$, with sensitivity $\mathscr{B}_1$, the utility loss $\sum_{i=1}^N f_i(\bar{\bm{x}}^K) - \sum_{i=1}^N f_i(\bm{x}^*)$ is logarithmically dependent on $d$ and upper bounded by $\widetilde{O}\big(\frac{\sqrt{M}NC_1\mathscr{B}_1}{m^{3/2}\epsilon}\big)$ for $(\epsilon, \delta)$-LDP.
	\end{thm}
	 As a final remark, it is noted that all the convergence analyses presented in this paper do not rely on any assumption on the constraint domain $\mathcal{C}$ but depend on the initialization $\bm{x}^0_i$, where in particular for the strongly convex case, the utility loss scales logarithmically with $f_i(\bm{x}^0_i)-f_i(\bm{x}^*)$. In addition, we present a privacy analysis of a variant of SGD: random coordinate descent (RCD) \cite{coordinate2015}. The difference between RCD and GD is that in each iteration, gradient descent is implemented only on one dimension randomly selected and thus complexity is reduced by $1/d$. In practice, many functions $f_i(\cdot)$ may have a much smaller coordinate smooth constant $\hat{M}$, which is defined as $\max_{l, \bm{x} \in \mathcal{C}} |\nabla f(\bm{x}+te_l)[l] - \nabla f(\bm{x})[l]| < \hat{M}|t|$ for $l \in [1:d]$ and $e_l$ is the $l^{th}$ standard basis of $\mathbb{R}^d$. It is not hard to observe that $1 \leq {M}/{\hat{M}} \leq d$. Details about private RCD can be found in Appendix K. 

	\section{Conclusion}
	In this paper, we investigate LDP in interactive decentralized optimization from both asymptotic and non-asymptotic viewpoints. A framework to incorporate randomized parameters is proposed, which can reduce the privacy loss by a constant while preserving the convergence rate. Simulation results are provided to validate the theory. Theoretically, we study asymptotic utility loss with sensitivity in $l_{\infty}$, $l_2$ and $l_1$ norms. Without resort to any sparse assumptions on the constraint $\mathcal{C}$, we find that even under LDP, utility loss may not explicitly depend on dimensionality $d$ when the objective functions are in good structures, analogous to existing study on DP in constrained optimization. 
	\newpage 
   \bibliographystyle{unsrt}
    \bibliography{ref}

	\newpage
	\appendix 
	
	
	\section{Private Decentralized Optimization with Secret Sharing}
	In this section, we provide the proof of Theorem 2.1. We first give the details of the algorithm construction, where the main idea is to securely compute the sum of exchanged updates via secret sharing. Without loss of generality, we take the ADMM algorithm defined in (4) as an example and the modified protocol is presented as below. An illustration is presented as Fig. (A). 
	\begin{algorithm}
		\label{ss-admm}
		\caption{$A_i\bm{x}^k_i$ Exchange in Secret-sharing based ADMM}
		\begin{algorithmic}
			\STATE {\bfseries Input:} $A_i\bm{x}^k_i$, $i=1,2,...,N$, $p \in \mathbb{Z}$
			\STATE {\bfseries Agents $i=1,2,...,N$ do in parallel: }
			\STATE $v_i$ randomly splits $A_i\bm{x}^k_i$ into $N$ shares, $\bm{s}^k_{ [1:N]}$, such that
			$$ A_i\bm{x}^k_i = \sum_{j=1}^{N} \bm{s}^k_{ij} \mod p. $$
			\STATE $v_i$ sends $\bm{s}^k_{ij}$ to node $v_j$ while keeping $\bm{s}^k_{ii}$ as a secret.
			\STATE {\bfseries Agents $i=1,2,...,N$ do in parallel: }
			\STATE $v_i$ sums up $\bm{s}^k_{[1:N] }$ received as
			$$ \hat{\bm{s}}^k_i = \sum_{j=1}^{N} \bm{s}^k_{ji} \mod p. $$
			\STATE $v_i$ broadcasts $\hat{\bm{s}}^k_i$
			\STATE Reconstruct $\sum_{i=1}^{N} A_i\bm{x}^k_i = \sum_{i=1}^{N} \hat{\bm{s}}^k_i  \mod p$
		\end{algorithmic}
	\end{algorithm}
	\begin{figure}[h!]
		\centering
		\includegraphics[width=1\linewidth]{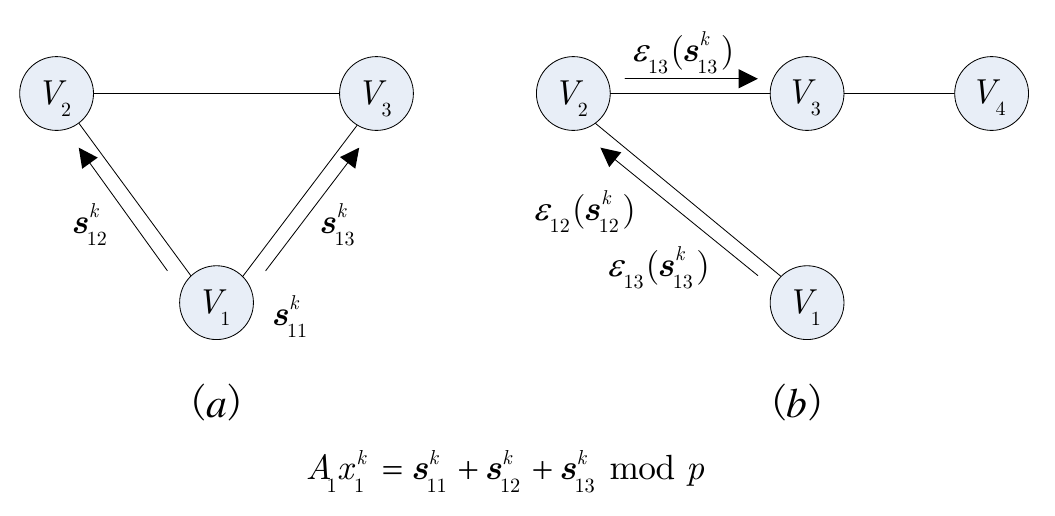}
		\captionsetup{labelformat=empty}
		\caption{Fig. (A) Secret-sharing based Exchange Protocol}
	\end{figure}
	In Algorithm 3, due to the linear constraint assumption, the procedure to update both $\bm{x}_i$ and $\bm{\lambda}$ only relies on the sum of $A_i\bm{x}^{k+1}_i$. Assume that in the $k $-th iteration, each participant has updated the states to $\bm{x}^{k}_i$. Let $p$ be a sufficiently large integer preselected such that $p > \Vert \sum_{i=1}^N A_i\bm{x}_i  \Vert_{\infty}$, i.e., $p$ is greater than the largest coordinate of $\sum_{i=1}^N A_i\bm{x}_i$ in absolute value. To share $A_i\bm{x}^{k}_i$ with neighbors, rather than exchanging directly, node $i$ randomly divides $A_i\bm{x}^{k}_i$ into $N$ shares, $\{\bm{s}^{k}_{i[1:N]}\}$, such that $\sum_{j=1}^{N} \bm{s}^{k}_{ij} = A_i\bm{x}^{k}_i$. Such division can be performed by randomly selecting $\bm{s}^{k}_{i[2:N]}$, and then $\bm{s}^{k}_{i1}$ is determined by $A_i\bm{x}^{k}_i-\sum_{j=2}^{N} \bm{s}^{k}_{ij}$.
	Then, $v_i$ sends the share $\bm{s^k_{ij}}$ to $v_j$, $j=[1:N] \backslash i$, while keeping $\bm{s}^k_{ii}$ to itself. After the exchange, each node $v_i$ still holds $N$ shares, $\bm{s}^k_{[1:N] i}$, of which one is from itself and the remaining are from the other $(N-1)$ neighbors. Then, each $v_i$ sums up all the shares held, denoted by $\hat{\bm{s}}^{k}_i = \sum_{j=1}^{N} \bm{s}^k_{ji}$ and broadcasts. Clearly, 
	\begin{equation}
	\sum_{i=1}^{N} \hat{\bm{s}}^k_i \equiv \sum_{i=1}^{N}  \sum_{j=1}^{N} \bm{s}^k_{ji} \equiv  \sum_{i=1}^{N}\sum_{j=1}^{N}   \bm{s}^k_{ij}  \equiv \sum_{i=1}^{N} A_i \bm{x}^{k}_i \mod p.
	\end{equation}
	Moreover, $\bm{x}^{k}_i$ can be reconstructed if and only if $\bm{s}^{k}_{i[1:N]}$ are all collected. For no more than $(N-2)$ colluding nodes, there always exists one share among $\bm{s}^{k}_{i[1:N]}$ which cannot be inferred by any $v_j$, $j \not =i$ and thus the scheme proposed is information theoretically secure. Especially, if we are concerned with a network adversary that can eavesdrop on all communication, we need to provide secure communication channels between each pair of nodes. A given pair of nodes $v_i$ and $v_j$ can encrypt shares $\bm{s}^k_{ij}$ and $\bm{s}^k_{ji}$ using a shared symmetric key prior to exchange. An example is given as Fig. A (b) (in a relay fashion), where we denote the ciphertext of $\bm{s}^k_{ij}$ by $\epsilon(\bm{s}^k_{ij})$. Relying on secret sharing, the proposed ADMM achieves privacy without any compromise in utility, while it comes with an additional round of data exchange in each iteration.
	
	A rigorous proof is given as follows. Considering $A_{i_0}\bm{x}^k_{i_0}$ for arbitrary $i_0 \in \{1,2,...,N\}$, based on the definition of secret splitting, one may reconstruct $A_{i_0}\bm{x}_{i_0}$ if and only all the $N$ shares have been collected (and decrypted properly in the encryption case). In the first step of Algorithm 3, where $(N-1)$ random shares have been distributed to the remaining $(N-1)$ nodes, there should exist at least one honest node, denoted by $v_{i_1}$ with shares $\bm{s}^k_{i_0i_1}$. Then, in the second step, each node sums up all the shares received as $\bm{\hat{s}}^k_{[1:N]}$ and broadcasts. It is clear that $$\bm{\hat{s}}^k_{i_0}= \sum_{j=1}^{N} \bm{s}^k_{ji_0} \mod p,$$ of which the reconstruction requires both $\bm{s}^k_{i_1i_0}$ and $\bm{s}^k_{i_0i_0}$, while $\bm{s}^k_{i_0i_0}$ is a secret of $v_{i_0}$ and $\bm{s}^k_{i_1i_0}$ is a secret between $v_{i_1}$ and $v_{i_0}$. With the assumption that $v_{i_1}$ is honest, for $v_i$, $i \not = i_1 , i_0$, from $\bm{\hat{s}}_{i_0}$, it is impossible to infer either $\bm{s}^k_{i_1i_0}$ or $\bm{s}^k_{i_0i_0}$. With a similar reasoning, since $N \geq 3$, the reconstruction of $\bm{\hat{s}}^k_{i_0}$ is also determined by some $\bm{s}^k_{ii_0}$ for $i \not = i_0, i_1$, which is unknown to $v_{i_1}$. Thus, $v_{i_1}$ cannot infer $\bm{s}^k_{i_0i_0}$ either. In a nutshell, either for $v_{i}$, $i \not =i_1, i_0$ or $v_{i_1}$, at least one share, i.e., $\bm{s}^k_{i_0i_0}$, cannot be inferred and thus $A_{i_0}\bm{x}^k_{i_0}$ is secure to at most $(N-2)$ colluding nodes. 
	
	\section{Proof of Lemma 2.1}
	It is noted that for an $\bm{X} =\bm{x}^{[0:K]}_{[1:N]}$ observed, since there is no prior on the inputs $\mathscr{D}=\{f_1, ... ,f_i, ... ,f_N\}$ and $\mathscr{D'}=\{f_1, ... ,\hat{f}_i, ... ,f_N\}$,
	\begin{equation}
	\label{distri-difference-bound-1}
	\frac{ P\big( \mathscr{D}  \big| \bm{x}^{[0:K]}_{[1:N]} \big)}{P\big( \mathscr{D'}  \big| \bm{x}^{[0:K]}_{[1:N]}  \big)} = \frac{ P\big( \bm{x}^{[0:K]}_{[1:N]}  \big| \mathscr{D}\big)}{P\big(  \bm{x}^{[0:K]}_{[1:N]}  \big|  \mathscr{D'}\big)} =  \frac{ P\big( \bm{x}^{0}_{[1:N]}  \big| \mathscr{D} \big)  \prod_{k=1}^{K} P\big( \bm{x}^{k}_{[1:N]}  \big| \mathscr{D},\bm{x}^{[0:k-1]}_{[1:N]}\big)}{  P\big( \bm{x}^{0}_{[1:N]}  \big| \mathscr{D'} \big) \prod_{k=1}^{K} P\big( \bm{x}^{k}_{[1:N]}  \big|  \mathscr{D'},\bm{x}^{[0:k-1]}_{[1:N]} \big)},
	\end{equation}
	It is noted that $\mathscr{D}$ and $\mathscr{D'}$ differ in $f_i$ and $\hat{f}_i$, to which the distribution of $\bm{x}_{{[1:N] \backslash i}}$ is invariant, and $\bm{x}^k_i$ only depends on the private function of agent $i$ and $\bm{x}^{[1:k-1]}_{[1:N]}$. On the other hand, the initialization of $\bm{x}^{0}_{[1:N]}$ is independent of the dataset. Thus, (\ref{distri-difference-bound-1}) can be further simplified as 
	$$ \prod_{k=0}^{K-1} \frac{ P\big( \bm{x}^{k+1}_{i}  \big| f_i, \bm{x}^{[0:k]}_{[1:N]} \big)} {P\big( \bm{x}^{k+1}_{i}  \big| \hat{f}_i, \bm{x}^{[0:k]}_{[1:N]}\big)} = \prod_{k=0}^{K-1} \prod_{l=1}^d \frac{ P\big( \bm{x}^{k+1}_{i}[l]  \big| f_i, \bm{x}^{[0:k]}_{[1:N]} \big)} {P\big( \bm{x}^{k+1}_{i}[l]  \big| \hat{f}_i, \bm{x}^{[0:k]}_{[1:N]}\big)},$$
	since the noise on each dimension is i.i.d. By taking the logarithm of the above equation and recalling the definition of $\epsilon(\bm{X})$ in Definition 1, the lemma follows.

	\section{Proof of Theorem 2.2 }
	When we assume $\bm{\Delta}^{k+1}_i$ is a Laplace distribution, the distributions $P( \bm{x}^{k+1}_{i}[j] | f_i, \bm{x}^{[0:k]}_{[1:N]})$ in either Algorithm 1 or 2 share a very similar structure. Both follow a mixture Laplace distribution with a random mean. In Algorithm 1 the mean is randomly distributed in an interval starting from $\bm{x}_i^k - D^{-1}_i\nabla f_i(\bm{x}^k_i) + D^{-1}_i\bm{\lambda}^{k}_i$ to $\frac{\sum_{ j \not = i}\bm{x}_j^k}{N-1} - D^{-1}_i\nabla f_i(\bm{x}^k_i) + D^{-1}_i\bm{\lambda}^{k}_i$, while in Algorithm 2 the mean is randomly distributed in an interval starting from $\frac{2}{N} \bm{x}^k_{i_k} + \frac{\sum_{j \not = i_k, \hat{i}_k} \bm{x}^k_j }{N} -\eta_{k+1}\nabla f_i(\frac{\sum_{i=1}^{N} \bm{x}^k_i}{N})$ to $\frac{2}{N} \bm{x}^k_{\hat{i}_k} + \frac{\sum_{j \not = i_k, \hat{i}_k} \bm{x}^k_j }{N} -\eta_{k+1}\nabla f_i(\frac{\sum_{i=1}^{N} \bm{x}^k_i}{N})$. Without loss of generality, we focus on Algorithm 1. Since $\bm{\Delta}^{k+1}_i$ on each dimension is i.i.d. in Lap$(0,\beta^{k+1})$, recalling Lemma 2.1, at iteration $ k+1 $, the bound on privacy loss in the $l$-th dimension can be expressed as
	\begin{equation}
	\label{average-dp-bound}
	\hat{\epsilon}^{k+1}_l (\bm{X}) = \sup_{\hat{f}_i \in \mathscr{F}_i } \bigg| \log \frac{P\big( \bm{x}^{k+1}_{i}[l]  \big| f_i, \bm{x}^{[0:k]}_{[1:N]} \big)} {P\big( \bm{x}^{k+1}_{i}[l]  \big| \hat{f}_i, \bm{x}^{[0:k]}_{[1:N]} \big)} \bigg| = \max_{|t| \leq D^{-1}_i\mathscr{B}_{\infty} } \bigg| \log  \frac{\int_{ \tau^{k+1}_i[l] } e^{-{\beta_{k+1}}{| \bm{x}^{k+1}_{i}[l]-Y | }} dY } {\int_{ \tau^{k+1}_i[l]+t } e^{-{\beta_{k+1}}{ | \bm{x}^{k+1}_{i}[l]-Y | }} dY }  \bigg|.
	\end{equation}
	We reformulate this problem as follows. For $X \in \mathbb{R}$, we consider
	\begin{equation}
	\label{difference-obj}
	\max_{|t| \leq D^{-1}_i\mathscr{B}_{\infty} } \bigg| \log \frac{\int_{0}^{\omega} {\beta_{k+1}} e^{-{\beta_{k+1}} {| X - Y| }} dY} {\int_{ t }^{ t +\omega} {\beta_{k+1}} e^{-{\beta_{k+1}}{|X -Y| }} dY } \bigg| ,
	\end{equation} 
	for some positive numbers $\omega, \beta_{k+1}$ and $\mathscr{B}_{\infty}$. Here, $\omega$ corresponds to the length of the interval.
	
	For a fixed $t$, $|t| \leq \mathscr{B}_{\infty}$, if $X \not  \in [0, \omega] \cup [t, \omega+t ]$, then
	$$   \bigg| \log \frac{\int_{0}^{\omega} {\beta_{k+1}} e^{-{\beta_{k+1}}{| X - Y| } dY}} {\int_{ t }^{ t +\omega} {\beta_{k+1}} e^{-{\beta_{k+1}}{|X -Y| } } dY } \bigg| =  \bigg| \log \frac{\int_{0}^{\omega } e^{-{\beta_{k+1}}{| X - Y| }} dY } {e^{{\beta_{k+1}}{t}} \int_{0}^{\omega} e^{-{\beta_{k+1}}{|X -Y| } } dY } \bigg| =\big|{\beta_{k+1}}{t}\big| \leq D^{-1}_i{\beta_{k+1}}{\mathscr{B}_{\infty}} .$$
	In the following, without loss of generality, we assume $X \in [0,\omega]$, then $\int_{0}^{\omega} {\beta_{k+1}} e^{-{\beta_{k+1}}{| X - Y| }} dY = 2-e^{-\beta_{k+1} X}-e^{-\beta_{k+1}(\omega-X)}$. First, supposing that $X \in [t, \omega+t]$, then $\int_{t}^{\omega+t} {\beta_{k+1}} e^{-{\beta_{k+1}}{| X - Y| }} dY = 2 - e^{-\beta_{k+1} (X-t)} - e^{-\beta_{k+1}(\omega+t-X)}$. To show $$ e^{-\beta_{k+1} |t| } \leq \frac{2 - e^{-\beta_{k+1} X} - e^{-\beta_{k+1}(\omega-X)}} {2 - e^{-\beta_{k+1} (X-t)} - e^{-\beta_{k+1}(\omega+t-X)}} \leq e^{\beta_{k+1} |t| }, $$ it is equivalent to showing
	\begin{equation}
	\label{both-in}
	\left\{ 
	\begin{aligned}
	& 2   e^{\beta_{k+1} |t| } - e^{-\beta_{k+1} X + \beta_{k+1} |t|} - e^{-\beta_{k+1}(\omega-X)+\beta_{k+1} |t|} \geq  2 - e^{-\beta_{k+1} (X-t)} - e^{-\beta_{k+1}(\omega+t-X)}, \\
	& 2 - e^{-\beta_{k+1} X} - e^{-\beta_{k+1}(\omega-X)} \leq 2e^{\beta_{k+1} |t| } - e^{-\beta_{k+1} (X-t) +\beta_{k+1} |t|} - e^{-\beta_{k+1}(\omega+t-X)+\beta_{k+1} |t|}.
	\end{aligned} 
	\right. 
	\end{equation}
	Due to the symmetry, we merely prove the case that when $t \geq 0$, where (\ref{both-in}) can be rewritten as, 
	\begin{equation}
	\label{both-in-2}
	\left\{ 
	\begin{aligned}
	& 2 e^{\beta_{k+1} t } - e^{-\beta_{k+1} (X - t)} - e^{-\beta_{k+1}(\omega-X-t)} \geq 2 - e^{-\beta_{k+1} (X-t)} - e^{-\beta_{k+1}(\omega+t-X)} ,\\
	& 2 - e^{-\beta_{k+1} X} - e^{-\beta_{k+1}(\omega-X)} \leq 2e^{\beta_{k+1} t } - e^{-\beta_{k+1} (X-2t)} - e^{-\beta_{k+1}(\omega-X)}.
	\end{aligned} 
	\right. 
	\end{equation}
	Clearly, for the first inequality, it suffices to show 
	\begin{equation}
	\label{2-1}
	2(e^{\beta_{k+1} t }-1) \geq (e^{2\beta_{k+1} t }-1) e^{-\beta_{k+1}(\omega+t-X)},
	\end{equation} 
	and it can be further simplified as $2e^{\beta_{k+1}(\omega+t- X)} \geq e^{\beta_{k+1} t}+1 $. Such a claim follows clearly as $  \omega - X \geq 0$ and $\alpha>0$. For the second inequality, with similar reasoning, it is equivalent to 
	\begin{equation}
	\label{2-2}
	2e^{\beta_{k+1} X} \geq e^{\beta_{k+1} t}+1,
	\end{equation} 
	which holds since $X \geq t$. At last, we consider $X \not \in [t, t+\omega]$. Still, due to the symmetry, we can assume $t > 0$ and $X < t$. Then, it is equivalent to show: 
	\begin{equation}
	\label{one-in}
	\left\{ 
	\begin{aligned}
	& 2   e^{\beta_{k+1} t } - e^{-\beta_{k+1} (X - t)} - e^{-\beta_{k+1}(\omega-X)+\beta_{k+1} t} \geq e^{-\beta_{k+1} (t-X)} - e^{-\beta_{k+1}(\omega+t-X)}, \\
	& 2 - e^{-\beta_{k+1} X} - e^{-\beta_{k+1}(\omega-X)} \leq e^{\beta_{k+1} X} - e^{-\beta_{k+1} (\omega-X) } .
	\end{aligned} 
	\right. 
	\end{equation}
	As for the first inequality, assume that $g(t) = 2 e^{\beta_{k+1} t } - e^{-\beta_{k+1} (X - t)} - e^{-\beta_{k+1} (t-X)}  -  e^{-\beta_{k+1}(\omega-X)+\beta_{k+1} t} + e^{-\beta_{k+1}(\omega+t-X)} $. It is noted that when $t=0$, $x$ should be also be 0 based on the assumption and $g(0) =0$. On the other hand,
	\begin{equation}
	\label{3-1}
	\frac{d g}{d t} = \beta_{k+1} \big(2 e^{\beta_{k+1} t } - e^{-\beta_{k+1} (X - t)} + e^{-\beta_{k+1} (t-X)} -  e^{-\beta_{k+1}(\omega-X)+\beta_{k+1} t} - e^{-\beta_{k+1}(\omega+t-X)}  \big).
	\end{equation}
	Since $X<\omega$, to show $g(t)$ is non-decreasing with respect to $t$, it suffices to show that,
	$$2 e^{\beta_{k+1} t } - e^{-\beta_{k+1} (X - t)} + e^{-\beta_{k+1} (t-X)} -  e^{-\beta_{k+1}(t-X)+\beta_{k+1} t} - e^{-\beta_{k+1}(t+t-X)} \geq 0.$$
	It is clear that $e^{\beta_{k+1} t } \geq e^{-\beta_{k+1} (X - t)}$ and $e^{-\beta_{k+1} (t-X)} \geq e^{-\beta_{k+1}(2t-X)}$ as both $X$ and $t$ are non-negative. Furthermore, $e^{\beta_{k+1} t } \geq e^{-\beta_{k+1}(t-X)+\beta_{k+1} t} = e^{\beta_{k+1} X}$ since $t \geq X$. Therefore, (\ref{3-1}) is non-negative. The second inequality of (\ref{one-in}) is exactly the AM-GM inequality that $$2  \leq e^{-\beta_{k+1} X}+ e^{\beta_{k+1} X}.$$  
	
	In a nutshell, we have proven that (\ref{difference-obj}) is upper bounded by $\max_{|t| \leq D^{-1}_i\mathscr{B}_{\infty}} |t\beta_{k+1}| = \beta_{k+1} D^{-1}_i \mathscr{B}_{\infty}$. Moreover, when $X$ belongs to the intersection of the two intervals, $(0,\omega)$ and $(t,\omega+t)$, the above inequalities are strict, i.e., (\ref{difference-obj}) is strictly smaller than $\beta_{k+1}D^{-1}_i \mathscr{B}_{\infty}$, which is the case if we fix all parameters to be constants. Similarly, by replacing $D^{-1}_i$ with $\eta_{k+1}$, we derive the proof for the case of Algorithm 2.

	\section{Proof of Theorem 2.3}
	We drop all the dependence on $i$, $k$ and $l$ for brevity. Following the normalization in the proof of Theorem 2.1, we still assume $x$ is a Laplace distribution of which the mean is uniformly distributed in $[0, \omega]$, conditional on all prior intermediate outputs. As the corollary of Theorem 2.1, 
	$$ \Theta(x) = \max_{|t| \leq \alpha\mathscr{B}_{\infty}} \bigg| \log  \frac{\int^{\omega}_{0} e^{-{\beta}{| x -y| }} dy } {\int^{\omega+t}_{t}  e^{-{\beta}{| x-y| }} dy }  \bigg|, $$
	where the maximization is achieved when $t$ either equals to $\alpha\mathscr{B}_{\infty}$ or $-\alpha\mathscr{B}_{\infty}$. To quantify $\gamma$, it suffices to calculate
	\begin{equation}
	\label{adv_1}
	\int_{-\infty}^{\infty} \int^{\omega}_{0} \Theta(x) \frac{\beta}{2\omega}e^{-{\beta}{| x-y| }} dy dx, 
	\end{equation}
	since the probability density function of $x$ is $\int_{0}^{\omega} \frac{\beta}{2\omega}e^{-\beta|x-y|}dy$. With the concavity of $\log(\cdot)$, (\ref{adv_1}) is upper bounded by 
	\begin{equation}
	\label{concave}
	\log \int_{-\infty}^{\infty} \int^{\omega}_{0} \max_{t = \pm \alpha\mathscr{B}_{\infty}} \bigg\{ \frac{\int^{\omega}_{0} e^{-{\beta}{| x -z| }} dz } {\int^{\omega+t}_{t}  e^{-{\beta}{| x-z| }} dz }, \frac{\int^{\omega+t}_{t}  e^{-{\beta}{| x-z| }}dz }{\int^{\omega}_{0} e^{-{\beta}{| x -z|  }} dz } \bigg\} \frac{\beta}{2\omega} e^{-{\beta}{| x-y| }} dx dy. 
	\end{equation}
	Now we take a closer look into $\max_{t = \pm \alpha\mathscr{B}_{\infty}} \big\{ \frac{\int^{\omega}_{0} e^{-{\beta}{| x -z| }} dz } {\int^{\omega+t}_{t}  e^{-{\beta}{| x-z| }} dz }, \frac{\int^{\omega+t}_{t}  e^{-{\beta}{| x-z| }}dz }{\int^{\omega}_{0} e^{-{\beta}{| x -z|  }}dz  } \big\}$. Still from the corollary of Theorem 2.1, once $\bm{x}^{k+1} \not \in [0, \omega]$, $\max_{t = \pm \alpha\mathscr{B}_{\infty}} \big\{ \frac{\int^{\omega}_{0} e^{-{\beta}{| x -z| }} dz } {\int^{\omega+t}_{t}  e^{-{\beta}{| x-z| }}dz  }, \frac{\int^{\omega+t}_{t}  e^{-{\beta}{| x-z| }}dz }{\int^{\omega}_{0} e^{-{\beta}{| x -z|  }}dz  } \big\} = e^{\alpha\beta\mathscr{B}_{\infty}}.$ 
	
	Since we assume $\omega > \alpha\mathscr{B}_{\infty}$,  it is not hard to observe that
	\begin{itemize}
		\item $x \in [0,\frac{\omega-\alpha\mathscr{B}_{\infty}}{2}]$, $\max_{t = \pm \alpha\mathscr{B}_{\infty}} \bigg\{ \frac{\int^{\omega}_{0} e^{-{\beta}{| x -z| }} dz } {\int^{\omega+t}_{t}  e^{-{\beta}{| x-z| }} dz}, \frac{\int^{\omega+t}_{t}  e^{-{\beta}{| x-z| }} dz}{\int^{\omega}_{0} e^{-{\beta}{| x -z|  }} dz } \bigg\} = \frac{\int^{\omega+t}_{t}  e^{-{\beta}{| x-z| }}dz }{\int^{\omega}_{0} e^{-{\beta}{| x -z|}} dz} \bigg|_{t=-\alpha\mathscr{B}_{\infty}};$
		\item $x \in [ \frac{\omega-\alpha\mathscr{B}_{\infty}}{2},\frac{\omega}{2}]$, $\max_{t = \pm \alpha\mathscr{B}_{\infty}} \bigg\{ \frac{\int^{\omega}_{0} e^{-{\beta}{| x -z| }} dz } {\int^{\omega+t}_{t}  e^{-{\beta}{| x-z| }} dz}, \frac{\int^{\omega+t}_{t}  e^{-{\beta}{| x-z| }} dz}{\int^{\omega}_{0} e^{-{\beta}{| x -z|  }} dz } \bigg\} = \frac{\int^{\omega }_{0}  e^{-{\beta}{| x-z| }}dz }{\int^{\omega+t}_{t} e^{-{\beta}{| x -z|}} dz} \bigg|_{t=-\alpha\mathscr{B}_{\infty}};$
		\item $x \in [\frac{\omega}{2},\frac{\omega+\alpha\mathscr{B}_{\infty}}{2}]$, $\max_{t = \pm \alpha\mathscr{B}_{\infty}} \bigg\{ \frac{\int^{\omega}_{0} e^{-{\beta}{| x -z| }} dz } {\int^{\omega+t}_{t}  e^{-{\beta}{| x-z| }} dz}, \frac{\int^{\omega+t}_{t}  e^{-{\beta}{| x-z| }} dz}{\int^{\omega}_{0} e^{-{\beta}{| x -z|  }} dz } \bigg\} = \frac{\int^{\omega}_{0}  e^{-{\beta}{| x-z| }}dz }{\int^{\omega+t}_{t} e^{-{\beta}{| x -z|}} dz} \bigg|_{t= \alpha\mathscr{B}_{\infty}};$
		\item $x \in [\frac{\omega+\alpha\mathscr{B}_{\infty}}{2},\omega]$, $\max_{t = \pm \alpha\mathscr{B}_{\infty}} \bigg\{ \frac{\int^{\omega}_{0} e^{-{\beta}{| x -z| }} dz } {\int^{\omega+t}_{t}  e^{-{\beta}{| x-z| }} dz }, \frac{\int^{\omega+t}_{t}  e^{-{\beta}{| x-z| }}dz }{\int^{\omega}_{0} e^{-{\beta}{| x -z| }}dz  } \bigg\} = \frac{\int^{\omega+t}_{t}  e^{-{\beta}{| x-z| }}dz }{\int^{\omega}_{0} e^{-{\beta}{| x -z| }}dz} \bigg|_{t=\alpha\mathscr{B}_{\infty}}.$
	\end{itemize}
	Thus, fortunately, we can avoid the complicated integral at least in $x \in [0,\frac{\omega}{2}-\alpha\mathscr{B}_{\infty}]$, or $x \in [\frac{\omega}{2}+\alpha\mathscr{B}_{\infty},\omega]$ where it is simplified to $O(\int^{\omega+t}_{t}  e^{-{\beta}{| x-y| }}dy)$. Now we can split $\mathbb{R}$ into three parts, $(-\infty,0) \cup (\omega,\infty)$, $[0, \frac{\omega-\alpha\mathscr{B}_{\infty}}{2}] \cup [ \frac{\omega+\alpha\mathscr{B}_{\infty}}{2}, \omega]$ and the rest $ (\frac{\omega-\alpha\mathscr{B}_{\infty}}{2},\frac{\omega+\alpha\mathscr{B}_{\infty}}{2}) $. To avoid the tedious term when $x \in  (\frac{\omega-\alpha\mathscr{B}_{\infty}}{2},\frac{\omega+\alpha\mathscr{B}_{\infty}}{2}) $, we simply use the global upper bound to simplify them to derive a closed-form expression but one may obtain the expression of $\gamma$ exactly with the same reasoning. Note the symmetry on $t=\pm \alpha\mathscr{B}_{\infty}$, (\ref{concave}) is upper bounded by
	\begin{equation}
	\log \bigg\{ e^{\omega\alpha\mathscr{B}_{\infty}} \bigg[1-2\int_{0}^{(\omega-\alpha\mathscr{B}_{\infty})/2} \int_{0}^{\omega} \frac{\beta}{2\omega}e^{-\beta|x-y|}dy dx\bigg] + 2 \int_{0}^{(\omega-\alpha\mathscr{B}_{\infty})/2} \int_{-\alpha\mathscr{B}_{\infty}}^{\omega-\alpha\mathscr{B}_{\infty}}\frac{\beta}{2\omega}e^{-{\beta}{| x -y| }} dydx \bigg\}.
	\end{equation}
	
	
	\section{Simulation Results}
	
	We test the proposed schemes and state-of-art approaches on regularized empirical risk minimization (ERM) tasks. We use the standard $Adult$ dataset from the UCI Machine Learning Repository. For simplicity, we call the task as UCI in the following. In UCI, the dataset consists of demographic records, including age, sex and income etc. in 15 total features. We try to predict whether the annual income of an individual is above $50k$. After processing of the data, we remove all individuals with missing values and normalize both columns (features) and rows (individuals) while converting labels $\{\geq 50k, <50k\}$ to $\{0,1\}$. The training samples are denoted by $\{\bm{z}^i_j \in \mathbbm{R}^{14}, \mathscr{L}^i_j \in\{0,1\} | i=1,\cdots,N, j =1,\cdots,b_i \}$. Consistent with \cite{ICML2018}, \cite{ADMM-DP}, we select $\mathscr{L}(x) = \log(1+\exp(-x))$. Thus, $N$ agents are collaboratively solving the following logistic regression: 
	$$ \min_{\bm{x}} \sum_{i=1}^{N} f_i(\bm{x}) = \sum_{i=1}^{N} \bigg(\frac{1}{b_i}\sum_{j=1}^{b_i} \log(1+\exp(-\mathscr{L}^i_j \bm{x}^T\bm{z}^i_{j}))+\frac{1}{2}\norm{\bm{x}}^2\bigg). $$
	UCI is run with different parameter settings. 10 independent runs of each algorithm for comparison are performed and each agent is randomly assigned $100$ samples from the dataset. In each run, the communication graph is randomly generated using the given $N$ and the number of edges $|\mathscr{E}|$. 
	
	In UCI, four examples (a), (b), (c) and (d) are provided. We uniformly assume that $D_i=D=10$ and $\zeta=0.5$ in all cases for Algorithm 1. In the case of private ADMM, previous works all assume fixed parameters in the optimization protocol. In \cite{ADMM-DP}, the Lagrangian multiplier at the beginning of each iteration is perturbed, while \cite{optimal2019} considers the output perturbation at the end of each iteration. Further, in \cite{ICML2018}, the authors introduce a sequence of increasing step penalty, which can bring better utility-privacy tradeoff empirically. 
	For \cite{ADMM-DP}, \cite{optimal2019} with constant fixed penalty, we assume $\Gamma_i = 0.5D$ and $\rho_i = \frac{0.5D}{|\mathscr{N}_i|}$, corresponding to the expectation of the penalty terms in Algorithm 1. Here $\mathscr{N}_i$ denotes the neighbors of agent $i$.
	As for \cite{ICML2018}, we follow their setting that $\Gamma^k_i = 0.5\times1.02^k|\mathscr{N}_i|$ and $\rho^k_i = 0.5\times1.02^k$. \footnote{We do not optimize the increasing penalty here but we find that in some cases by proper selection, a privacy loss reduction can be achieved empirically at a cost of relatively small utility compromise. Such techniques can also be applied in our algorithms.}
	
	In the privacy part, with the same assumption in \cite{ICML2018}, we assume $f_i$ and $\hat{f}_i$ may only differ in one sample and thus, due the normalization, ${\mathscr{B}_{\infty}} =\frac{1}{b_i}=0.01$, and $\bm{J} = \frac{2.8}{DB_i}$, the Jacobian constant required by \cite{ICML2018} in their privacy analysis. It is noted that derivative of $\mathscr{L}$ is within $(-1,0]$, while the privacy analysis of \cite{optimal2019} requires a global sensitivity on that of $\nabla \mathscr{L}$. This makes their bound in this example too loose and we omit their privacy loss bound in our simulation. Following the setting of \cite{ICML2018}, \cite{optimal2019}, we also use a diminishing noise by selecting $\beta_k = 1.02^k$. The results of Example (a) are illustrated in Fig. 2, where $N=10$, $|\mathscr{E}|=20$. The accuracy logarithm defined by $ \log \norm{(\bm{x}^k_i-\bm{x}^*_i)/d}$, across 100 iterations averaged across 10 runs. The difference between the best and the worst accuracy over 100 runs is also marked. In example (b), with illustration shown in Fig. (B), under the same setting, we test algorithms in a large-scale case where $N=100$, $|\mathscr{E}|=200$. 
	\begin{figure}[h!]
		\centering
		\includegraphics[width=1\linewidth]{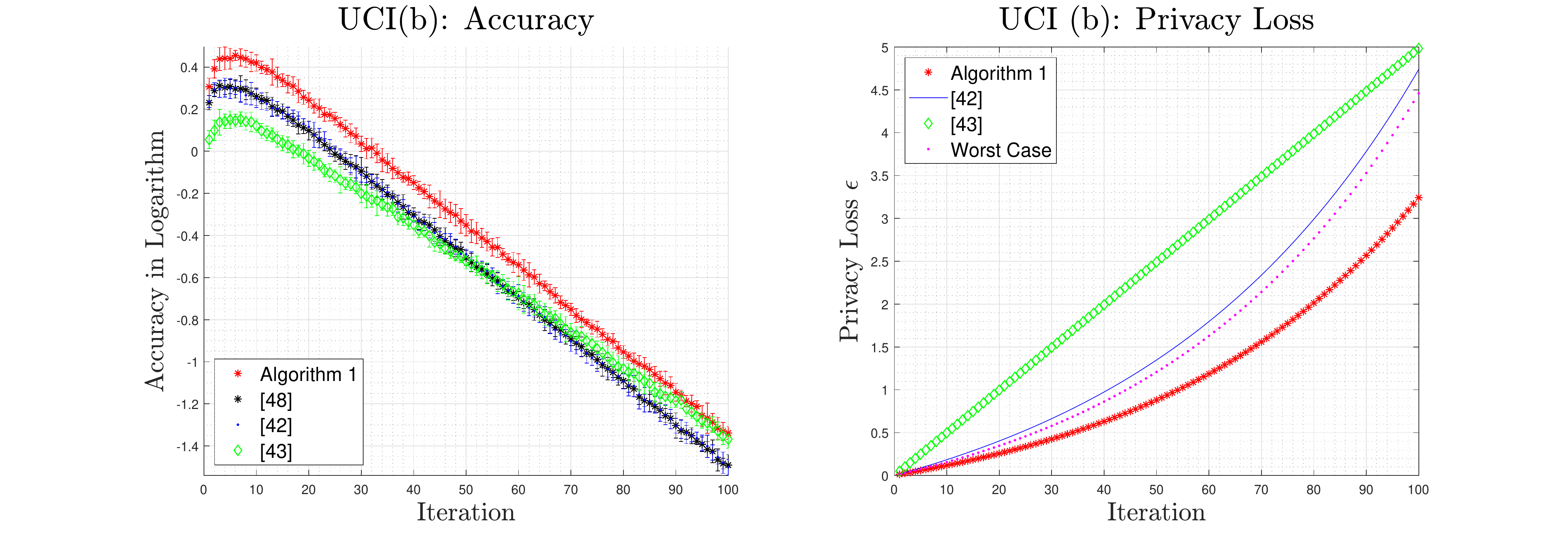}
		\captionsetup{labelformat=empty}
		\caption{Fig. (B) $N=100$, $\mathscr{E}=200$}
	\end{figure}

	In Example (c), we present an interesting variant to Algorithm 2. With Theorem 2.3, it is clear that under our framework, a larger interval can produce a better privacy amplification. Instead of the construction in Algorithm 2, one can construct a random aggregator by more aggressively utilizing the divergence among $\bm{x}^k_{[1:N]}$. For an instance, let $\widetilde{\mathscr{N}}_i$ denote the neighbors of node $i$ including $i$. To aggregate $\bm{x}^{k}_{j \in \mathscr{N}_i}$, for each dimension $l \in [1:d]$, let $$\alpha^k_{i,\min} = \min_{j \in \widetilde{\mathscr{N}}_i} \bm{x}^{k}_{j}[l],~~\text{and}~~\alpha^k_{i,\max} = \max_{j \in \widetilde{\mathscr{N}}_i} \bm{x}^{k}_{j}[l].$$ We consider the updating subroutine (\ref{GD_new_1}) in Algorithm 2*, i.e., each coordinate of $\bm{x}^{k+1}_i[l]$ is uniformly selected between the largest $\max_{j \in \widetilde{\mathscr{N}}_i} \bm{x}^{k}_{j}[l]$ and the smallest $\min_{j \in \widetilde{\mathscr{N}}_i} \bm{x}^{k}_{j}[l]$. 
	\begin{algorithm}[t]
		\caption*{Algorithm 2*}
		\begin{algorithmic}
			\STATE {\bfseries Input:} Local functions $f_i$ and a diminishing sequence $\{\eta_{k}\}$
			\STATE Initialize $\bm{x}^0_{[1:N]}$. 
			\FOR{$k=0,1,2, ...,K-1 $}
			\STATE {\bfseries Agents $i=1$ to $N$} do in parallel :\\
			\FOR{$l=1,2,...,d$}
			\STATE Randomly and independently generating a weight $w$ within $(0,1)$ and then updating $\bm{x}^{k+1}_i[l]$:
			\begin{equation}
			\label{GD_new_1}
			\bm{x}^{k+1}_i[l] := w \cdot \alpha^k_{i,\min} +(1-w) \cdot \alpha^k_{i,\max} - \eta_{k+1}\nabla f_i(\bm{x}^k_i) +\bm{\Delta}^{k+1}_i[l].
			\end{equation}
			\STATE Exchange $\bm{x}^{k+1}_i $ with Neighbors 
			\ENDFOR
			\ENDFOR
		\end{algorithmic}
	\end{algorithm}
	
	Different from Algorithm 1 and 2 which are controllable, aggregation in Algorithm 2* will bring some compromise in convergence. Following \cite{SGD2015}, we also select a diminishing step size, where for \cite{SGD2015}, $\eta_k = 0.9^k$ and in Algorithm 2*, $\eta_k = 0.93^k$ for balance. Not surprisingly, Algorithm 2* has a worse convergence at the beginning since information from neighbors is less efficiently merged but finally \cite{SGD2015} and Algorithm 2* achieve almost the same utility loss. However, the privacy loss of Algorithm 2* is only $30\%$ of \cite{SGD2015}, as shown in Fig (C). Since Algorithm 1 applies a constant step size, it has better accuracy but worse privacy loss, where all noises are fixed to be $\beta_k=1.02^k$. 
	\begin{figure}[h!]
		\centering
		\includegraphics[width=1\linewidth]{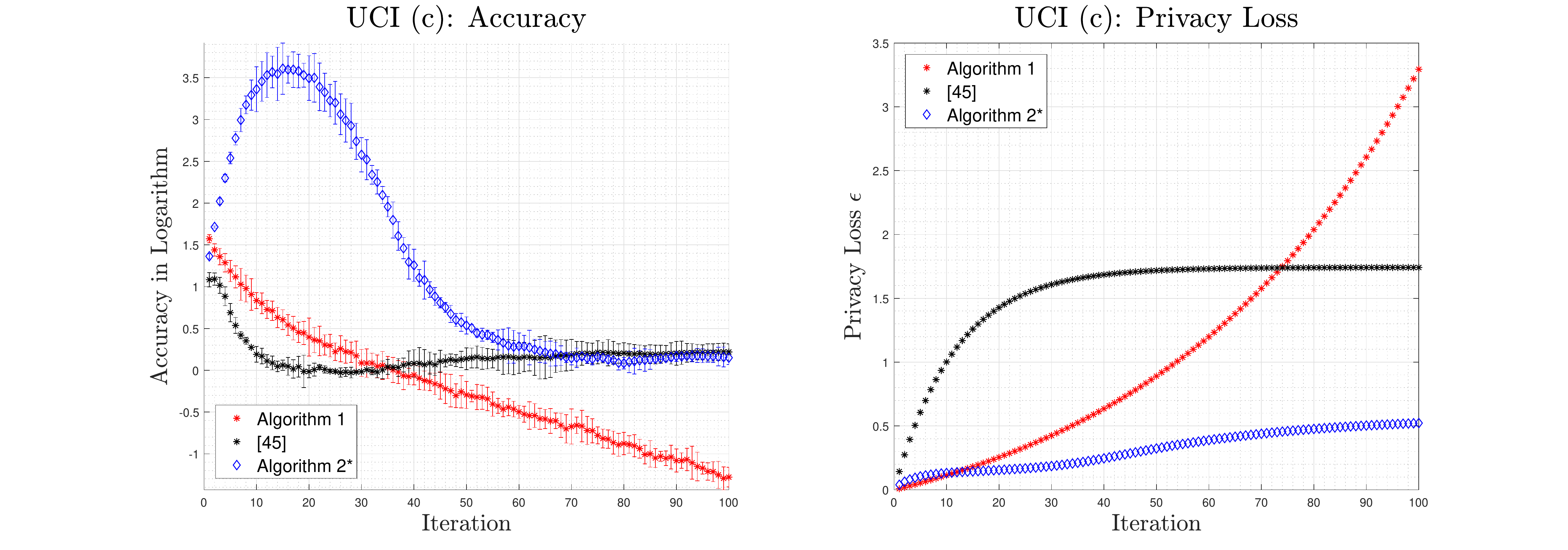}
		\captionsetup{labelformat=empty}
		\caption{Fig. (C) $N=10$, $\mathscr{E}=40$}
	\end{figure}
	
	Finally, we provide the performance of non-private optimization in UCI. With the same parameter setting as before, we test Algorithm 1 without first-order approximation and conventional ADMM with fixed parameters. In addition, we set the step size $\eta_k = 0.95^k$ and test Algorithm 2 and conventional decentralized GD (DGD) with fixed parameters. The performance is illustrated in Fig. D. In Fig. (D-1), the graph is randomly generated with $N=10$ and $|\mathscr{E}|=20$, same as Example (a) in Fig. 2; while in Fig. (D-2), $N=100$ and $|\mathscr{E}|=200$, same as Example (b). 
	\begin{figure}[htbp]
		\centering
		\begin{minipage}[t]{0.48\textwidth}
			\centering
			\includegraphics[width=6cm]{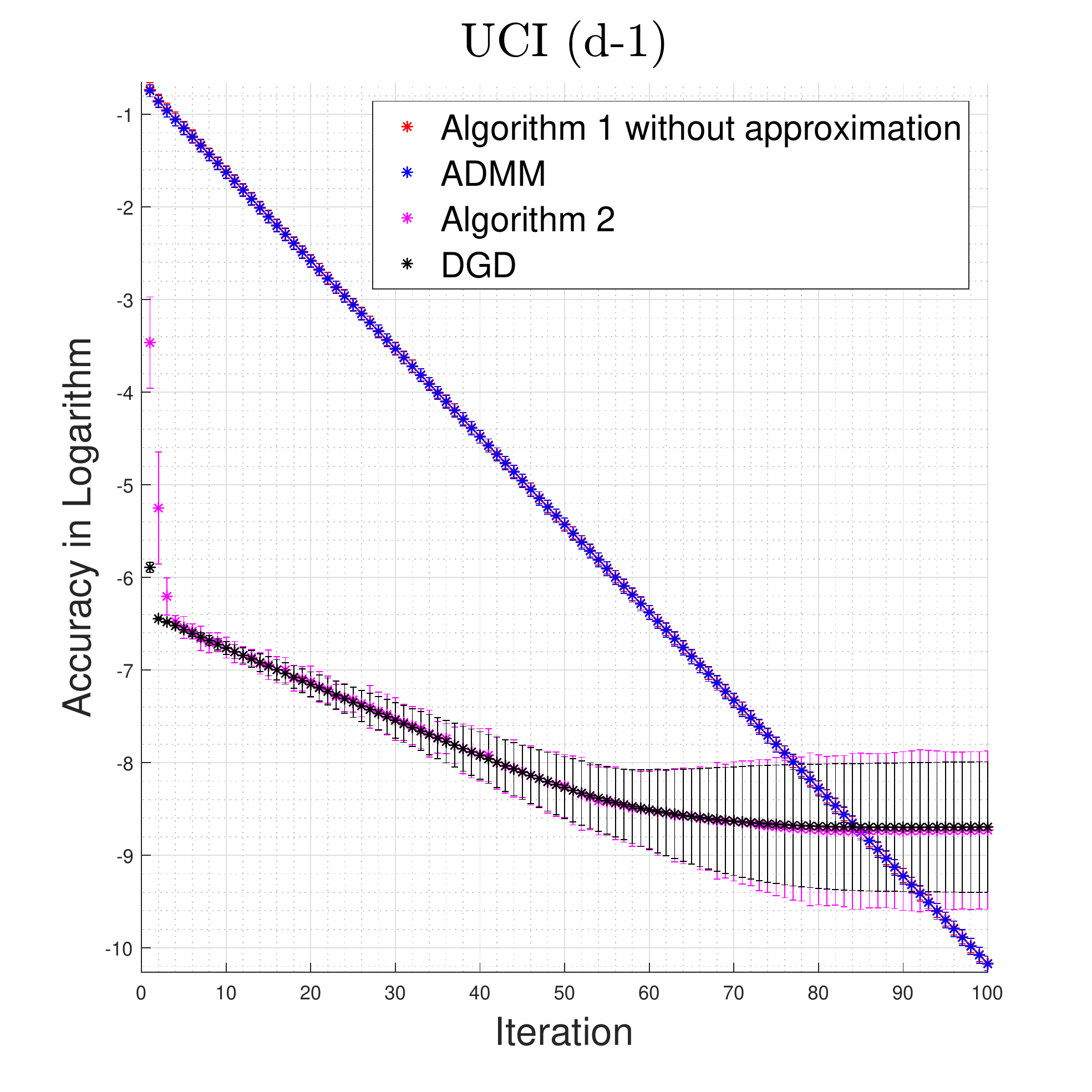}
			\captionsetup{labelformat=empty}
			\caption{Fig. (D-1) $N=10$, $|\mathscr{E}|=20$}
		\end{minipage}
		\begin{minipage}[t]{0.48\textwidth}
			\centering
			\includegraphics[width=6cm]{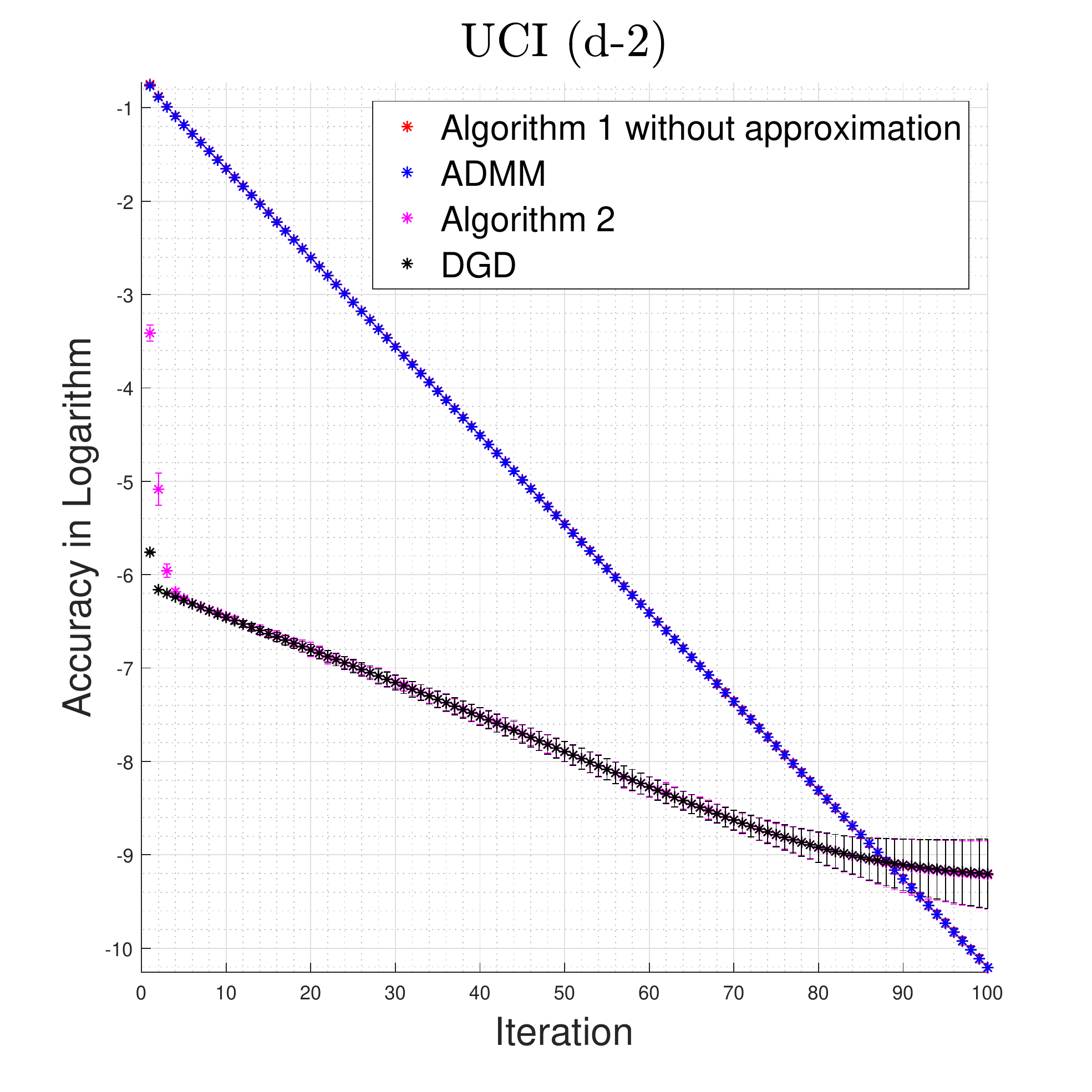}
			\captionsetup{labelformat=empty}
			\caption{Fig. (D-2) $N=100$, $|\mathscr{E}|=200$}
		\end{minipage}
	\end{figure}
	
	From Fig. D, it is clear, randomization defined in Algorithm 1 and 2 does not incur accuracy loss, which is consistent with our analysis.

	\section{Proof of Theorem \ref{varying-rho}}
	Since the proximal term $\norm{\bm{x}_i-\bm{x}_i^k}_{\bm{\Gamma}_i^{k+1}}^2$ is required to be nonnegative, the matrix $\bm{\Gamma}_i^{k+1}$ should be positive definite. With $D_i \cdot \bm{I} = A^T_i\bm{\rho}^{k+1}_iA_i + \bm{\Gamma}^{k+1}_i$, we just need to guarantee that the $D_i$ satisfy $D_i - \sigma_{\max}(A_i^T\bm{\rho}_i^{k+1}A_i)>0$ where $\sigma_{\max}(Z)$ and $\sigma_{\min}(Z)$ denote the maximal and the minimal non-zero singular value of Z, respectively. It leads to $D_i>\rho_{i,\max}^{k+1}\sigma^2_{i,\max}$ where $\sigma_{i,\max}$ is the largest singular value of $A_i$ and $\rho_{i,\max}^{ k+1 }$ is the maximum diagonal element of $\bm{\rho}_i^{k+1}$.
	
	To show the linear convergence, it suffices to determine $p >0$ such that,
	\begin{equation}
	\norm{\bm{u}^k - \bm{u}^*}^2 \geq (1+p ) \norm{\bm{u}^{k+1}-\bm{u}^*} ^2,
	\end{equation}
	which can be reformulated as 
	\begin{equation}
	\norm{\bm{u}^{k} - \bm{u}^{*}}^2 -\norm{\bm{u}^{k+1}-\bm{u}^{*}}^2 \geq p \norm{\bm{u}^{k+1} - \bm{u}^{*}}^2. 
	\end{equation}
	With the strong convexity, 
	\begin{equation}
	\label{basis}
	\langle \bm{x}-\bm{y}, \nabla f_i(\bm{x}) - \nabla f_i(\bm{y}) \rangle \geq m_i \norm{\bm{x}-\bm{y}}^2.
	\end{equation}
	And from (\ref{time-varying-x}), we have 
	\begin{equation}
	\label{gradient-1}
	\nabla f_i(\bm{x}^{k+1}_i) =  A^T_i(\bm{\lambda}^{k} - \bm{\rho}^{k+1}_i(A_i\bm{x}^{k+1}_i+ \sum_{j\not =i} A_j \bm{x}^{k}_j-\bm{c}) )+ \bm{\Gamma}_i^{k+1}(\bm{x}^{k}_i-\bm{x}^{k+1}_i).
	\end{equation}
	Also from the KKT condition, for the optimal states $\bm{\lambda}^*$ and $\bm{x}^* = (\bm{x}^*_1, ...,\bm{x}^*_N)$
	\begin{equation}
	\nabla f_i(\bm{x}^{*}_i) =  A^T_i\bm{\lambda}^*, \sum_{i=1}^{N} A_i\bm{x}^*_i = \bm{c}.
	\end{equation}
	Substituting the above equations into (\ref{basis})
	\begin{equation}
	(\bm{x}^{k+1}_i-\bm{x}^*_i)^T( A^T_i(\bm{\lambda}^{k} - \bm{\lambda}^{*}) - A_i^T\bm{\rho}^{k+1}_i(A_i(\bm{x}^{k+1}_i-\bm{x}^{k}_i)+ \sum_{ j=1}^{N} A_j( \bm{x}^{k}_j-\bm{x}^{*}_j)) + \bm{\Gamma}_i^{k+1}(\bm{x}^{k}_i-\bm{x}^{k+1}_i))\geq m_i \norm{\bm{x}^{k+1}_i-\bm{x}^{*}_i}^2.
	\end{equation}
	
	Summing up all for each $i$, it is noted that $\sum_{i=1}^{N} A_i(\bm{x}^{k+1}_i-\bm{x}^*_i) = \frac{1}{ \zeta }(\bm{\lambda}^{k}-\bm{\lambda}^{k+1})$ and
	\begin{equation}
	\label{G-expand}
	\begin{aligned}
	&(\bm{u}^{k+1}- \bm{u}^*)^TG(\bm{u}^{k}- \bm{u}^{k+1})=  \frac{1}{\zeta}(\bm{\lambda}^{k+1}-\bm{\lambda}^*)^T(\bm{\lambda}^{k}-\bm{\lambda}^{k+1}) + \sum_{i=1}^{N} (\bm{x}^{k+1}_i- \bm{x}^*_i)^T(A^T_i\bm{\rho}^{k+1}_iA_i+\bm{\Gamma}^{k+1}_i) (\bm{x}^{k}_i- \bm{x}^{k+1}_i)   \\
	& \geq -\frac{1}{\zeta} \norm{\bm{\lambda}^k - \bm{\lambda}^{k+1}}^2  + \sum_{i=1}^{N} m_i\norm{\bm{x}^{k+1}_i - \bm{x}^{*}_i}^2 + (\sum_{i=1}^{N} \bm{\rho}^{k+1}_iA_i(\bm{x}^{k+1}_i- \bm{x}^*_i))^T  (\sum_{j=1}^N A_j(\bm{x}^k_j- \bm{x}^*_j)).
	\end{aligned}
	\end{equation}
	Here, let the matrix $G= diag(\{\bm{D}_1,... ,\bm{D}_N, \frac{1}{\zeta}\})$, where $\bm{D}_i = D_i\cdot \bm{I}$, then it suffices to show $\norm{\bm{u}^k - \bm{u}^*}^2_G - \norm{\bm{u}^{k+1}-\bm{u}^*}^2_G \geq  p  \norm{\bm{u}^{k+1}-\bm{u}^*}^2_G$. On the other hand, $\norm{\bm{u}^{k} - \bm{u}^*}^2_G - \norm{\bm{u}^{k+1}-\bm{u}^*}^2_G = 2(\bm{u}^{k+1} - \bm{u}^*)^TG(\bm{u}^{k} - \bm{u}^{k+1}) + \norm{\bm{u}^k- \bm{u}^{k+1}}^2_G$. Referring to (\ref{G-expand}), it is equivalent to figure out $ p$ such that, 
	\begin{equation}
	\label{linear-1}
	\begin{aligned}
	& 2 \sum_{i=1}^{N} m_i\norm{\bm{x}^{k+1}_i-\bm{x}^{*}_i}^2+  2(\sum_{i=1}^{N} \bm{\rho}^{k+1}_iA_i(\bm{x}^{k+1}_i- \bm{x}^{*}_i))^T (\sum_{i=1}^{N} A_i(\bm{x}^{k}_i- \bm{x}^{*}_i)) + \sum_{i=1}^N D_i\norm{\bm{x}^{k+1}_i-\bm{x}^{k}_i}^2 - \frac{1}{\zeta}\norm{\bm{\lambda}^{k+1}-\bm{\lambda}^k}^2  \\
	& \geq  p  (\sum_{i=1}^N D_i\norm{\bm{x}^{k+1}_i-\bm{x}^{*}_i}^2 +\frac{1}{ \zeta }\norm{\bm{\lambda}^{k+1}-\bm{\lambda}^{*}}^2).
	\end{aligned}
	\end{equation}
	From (\ref{gradient-1}) and with the fact that $\bm{x}^{k}_i- \bm{x}^{*}_i = \bm{x}^{k}_i-\bm{x}^{k+1}_i+\bm{x}^{k+1}_i- \bm{x}^{*}_i$, we get
	\begin{equation}
	\label{upper-linear}
	\begin{aligned}
	\norm{\bm{\lambda}^{k+1}-\bm{\lambda}^*}^2 \leq& \frac{1}{\sigma^2_{i,\min}} \norm{A_i^T (\bm{\lambda}^{k+1}-\bm{\lambda}^*) }^2 \\
	=& \frac{1}{\sigma^2_{i,\min} } \bigg\|\nabla f_i(\bm{x}^{k+1}_i) - \nabla f_i(\bm{x}^{*}_i) - A_i^T(\bm{\lambda}^{k}-\bm{\lambda}^{k+1}) - \bm{D}_i(\bm{x}^{k}_i-\bm{x}^{k+1}_i) + A_i^T\bm{\rho}^{k+1}_i\sum_{j=1}^{N}A_j( \bm{x}^{k}_j-\bm{x}^{*}_j)\bigg\| ^2 \\
	\leq&  \frac{5}{\sigma^2_{i,\min} } ( M_i\norm{\bm{x}^{k+1}_i-\bm{x}^*_i }^2 + \sigma^2_{i,\max} \norm{\bm{\lambda}^{k}-\bm{\lambda}^{k+1}}^2 + D^2_{i} \norm{ 
		\bm{x}^{k}_{i}-\bm{x}^{k+1}_i}^2 +  \\
	& \rho_{i,\max}^{2(k+1)}\sigma^2_{i,\max}\bigg\|\sum_{j=1}^{N}A_j( \bm{x}^{k}_j-\bm{x}^{k+1}_j)\bigg\|^2+\rho_{i,\max}^{2(k+1)}\sigma^2_{i,\max}\big\| \sum_{j=1}^{N}A_j( \bm{x}^{k+1}_j-\bm{x}^{*}_j)\big\|^2),
	\end{aligned}
	\end{equation}
	where $\sigma_{i,\min} $ is the smallest nonzero singular value of $A_i$. For simplicity, $\rho_{i,\max}^{2(k+1)}=(\rho_{i,\max}^{k+1})^2$.
	Now, we substitute (\ref{upper-linear}) to (\ref{linear-1}), and then it can be reformulated as 
	\begin{equation}
	\begin{aligned}
	& \sum_{i} (2m_i-\frac{5M_ip}{\zeta N\sigma^2_{i,\min} }-D_ip)\norm{\bm{x}^{k+1}_i-\bm{x}^{*}_i}^2 + \sum_{i} (D_i-\frac{5p D^2_i}{\zeta N\sigma^2_{i,\min}}) \norm{\bm{x}^{k}_i-\bm{x}^{k+1}_i}^2\\
	&- (\frac{1}{\zeta}+\frac{5p }{\zeta N}\sum_{i=1}^{N} \frac{\sigma^2_{i,\max}}{\sigma^2_{i,\min}})\norm{\bm{\lambda}^k-\bm{\lambda}^{k+1}}^2+ 2(\sum_{i=1}^{N} \bm{\rho}^{k+1}_iA_i(\bm{x}^{k+1}_i- \bm{x}^{*}_i))^T (\sum_{j=1}^{N} A_j(\bm{x}^{k}_j- \bm{x}^*_j)) \\
	& - \frac{5p}{\zeta N}\sum_{i=1}^{N}\frac{\rho_{i,\max}^{2(k+1)}\sigma^2_{i,\max}}{\sigma^2_{i,\min}}\big(   \bigg\| \sum_{j=1}^{N}A_j( \bm{x}^{k}_j-\bm{x}^{k+1}_j)\bigg\|^2 + \frac{1}{\zeta^2}\norm{ \bm{\lambda}^k-\bm{\lambda}^{k+1}}^2\big)  \geq 0.
	\end{aligned}
	\end{equation}
	Moreover, $2(\sum_{i=1}^{N} \bm{\rho}^{k+1}_iA_i(\bm{x}^{k+1}_i- \bm{x}^{*}_i))^T (\sum_{j=1}^{N} A_j(\bm{x}^{k}_j- \bm{x}^*_j))$ can be rewritten as 
	\begin{equation}
	\begin{aligned}
	&2(\sum_{i=1}^{N} \bm{\rho}^{k+1}_iA_i(\bm{x}^{k+1}_i - \bm{x}^{*}_i))^T (\sum_{j=1}^{N} A_j(\bm{x}^{k}_j - \bm{x}^{*}_j))\\
	=&2(\sum_{i=1}^{N} \bm{\rho}^{k+1}_iA_i(\bm{x}^{k+1}_i - \bm{x}^{*}_i))^T (\sum_{j=1}^{N} A_j(\bm{x}^{k}_j - \bm{x}^{k+1}_j)) +  \frac{2}{\zeta^2} (\bm{\lambda}^{k}-\bm{\lambda}^{k+1})^T\bm{\rho}^{0}(\bm{\lambda}^{k}-\bm{\lambda}^{k+1}) + \\&\frac{2}{\zeta}(\sum_{i=1}^{N} (\bm{\rho}^{k+1}_i-\bm{\rho}^{0})A_i(\bm{x}^{k+1}_i- \bm{x}^{*}_i))^T (\bm{\lambda}^{k}-\bm{\lambda}^{k+1})\\
	\geq& \frac{2\rho^0}{\zeta^2}\norm{\bm{\lambda}^{k}-\bm{\lambda}^{k+1}}^2 - \sum_{i=1}^{N}\alpha N\rho_{i,\max}^{2(k+1)}\sigma^2_{i,\max}  \norm{\bm{x}^{k+1}_i- \bm{x}^{*}_i}^2-\sum_{i=1}^{N} \frac{N\sigma^2_{i,\max} }{\alpha }\norm{\bm{x}^{k}_i - \bm{x}^{k+1}_i}^2 -\\
	&\sum_{i=1}^{N} \alpha  \breve{\rho}_{i,\max}^{2(k+1)}\sigma^2_{i,\max}  \norm{\bm{x}^{k+1}_i- \bm{x}^{*}_i}^2-\frac{N}{\alpha \zeta^2}\norm{\bm{\lambda}^{k}-\bm{\lambda}^{k+1}}^2,
	\end{aligned}
	\end{equation}
	where $\bm{\rho}^0 = \rho^0\cdot \bm{I}$ and $\breve{\rho}_{i,\max}^{ k+1} $ is the maximum diagonal element of matrix $\bm{\rho}_i^{k+1}-\bm{\rho}^0$.
	Further, we have the following AM-GM inequality
	\begin{equation}
	\begin{aligned}
	&\bigg\| \sum_{j=1}^{N}A_j( \bm{x}^{k}_j-\bm{x}^{k+1}_j)\bigg\|^2 
	\leq  N\sum_{j=1}^{N}\sigma^2_{j,\max} \big\|\bm{x}^{k}_j-\bm{x}^{k+1}_j\big\|^2.
	\end{aligned}
	\end{equation}
	Combining (29), (30) and (31), we find that it suffices to find out $p$ such that 
	\begin{equation}
	\begin{aligned}
	&\sum_{i=1}^{N}\bigg(2m_i-\frac{5M_ip}{\zeta N\sigma^2_{i,\min} }-D_ip - \alpha  N\rho_{i,\max}^{2(k+1)}\sigma^2_{i,\max} -  \alpha  \breve{\rho}_{i,\max}^{2(k+1)}\sigma^2_{i,\max}  \bigg) \norm{\bm{x}^{k+1}_i-\bm{x}^{*}_i}^2 + \\&\sum_{i=1}^{N} \bigg(D_i- \frac{N\sigma^2_{i,\max} }{\alpha }-\frac{5p D^2_i}{\zeta N\sigma^2_{i,\min}}  - \frac{5p\sigma^2_{i,\max}}{\zeta}\sum_{ j=1}^{N}\frac{\rho_{j,\max}^{2(k+1)}\sigma^2_{j,\max}}{\sigma^2_{j,\min}}\bigg)\norm{\bm{x}^{k }_i-\bm{x}^{k+1}_i}^2+\\
	& \bigg(\frac{2\rho^0}{\zeta^2} -\frac{N}{\alpha \zeta^2}-\frac{1}{\zeta}  -\frac{5p }{\zeta  }\sum_{i=1}^{N}\bigg(\frac{\rho_{i,\max}^{2(k+1)}}{\zeta^2}+\frac{1}{N}\bigg) \frac{\sigma^2_{i,\max}}{\sigma^2_{i,\min}} \bigg)\norm{\bm{\lambda}^{k}-\bm{\lambda}^{k+1}}^2 \geq 0.
	\end{aligned}
	\end{equation}
	Therefore, $p$ can be selected as 
	\begin{equation}
	\label{delta-selection}
	\min \bigg\{ \frac{ 2m_i- \alpha  N\rho_{i,\max}^{2(k+1)}\sigma^2_{i,\max} - \alpha \breve{\rho}_{i,\max}^{2(k+1)}\sigma^2_{i,\max} }{\frac{5M_i}{\zeta N\sigma^2_{i,\min}}+D_i },\frac{D_i-\frac{N\sigma^2_{i,\max}}{\alpha }}{\frac{5  D^2_i}{\zeta N\sigma^2_{i,\min}} + \frac{5 \sigma^2_{i,\max}}{\zeta}\sum_{ j=1}^{N}\frac{\rho_{j,\max}^{2(k+1)}\sigma^2_{j,\max}}{\sigma^2_{j,\min}}} , \frac{\frac{2\rho^0}{\zeta }-\frac{N}{\alpha\zeta  }-1}{ 5  \sum_{i=1}^{N}(\frac{\rho_{i,\max}^{2(k+1)}}{\zeta^2}+\frac{1}{N}) \frac{\sigma^2_{i,\max}}{\sigma^2_{i,\min}}  } \bigg\}.
	\end{equation}
	To guarantee that $p>0$, the parameters $\alpha$, $D_i$, $\rho^0$ and $\zeta$ should satisfy:
	\begin{equation}
	\label{admissable-range}
	\left\{ 
	\begin{aligned}
	&\alpha <\frac{2m_i}{N\rho_{i,\max}^{2(k+1)}\sigma^2_{i,\max}+ \breve{\rho}_{i,\max}^{2(k+1)}\sigma^2_{i,\max}},    \\
	&D_i >\max\{ \rho_{i,\max}^{k+1}\sigma^2_{i,\max},\frac{N\sigma^2_{i,\max}}{\alpha}  \}, \\
	&\rho^0 >\frac{N}{2\alpha},\\
	&\zeta < 2\rho^0 - \frac{N}{\alpha}.
	\end{aligned}
	\right.
	\end{equation}

	\section{Proof of Theorem \ref{first-order-converge}}
	
	Under the strong continuity of both $f_i$ and its gradient $\nabla f_i$, for any $\bm{x}$ and $\bm{y}$, $$ \norm{\nabla f_i(\bm{x})  - \nabla f_i(\bm{y})}^2 \leq M_i \norm{\bm{x}-\bm{y}}^2 ,$$ and we use the following fact, for any $\bm{z}$
	$$ f_i(\bm{x})- f_i(\bm{y})  \leq \nabla f_i^T(\bm{z}) (\bm{x}- \bm{y}) + \frac{M_i}{2} \norm{\bm{x} - \bm{z}}^2,$$ and with strong convexity we have
	\begin{equation}
	\begin{aligned}
	\frac{m_i}{2}\norm{\bm{x}^{k+1}_i-\bm{x}_i^*}^2+ \nabla f_i(\bm{x}_i^*)^T(\bm{x}^{k+1}_i-\bm{x}_i^*)    \leq f_i(\bm{x}^{k+1}_i)-f_i(\bm{x}_i^*)  \leq \nabla f_i^T(\bm{x}^k_i) (\bm{x}^{k+1}_i-\bm{x}_i^*) + \frac{M_i}{2}  \norm{\bm{x}^{k}_i-\bm{x}^{k+1}_i}^2. 
	\end{aligned}
	\end{equation}
	On the other hand, since $A^T_i\bm{\lambda}^* = \nabla f_i(\bm{x}_i^*)$, thus
	$$  \frac{m_i}{2} \norm{\bm{x}^{k+1}_i-\bm{x}_i^*}^2 \leq (\bm{x}^{k+1}_i-\bm{x}_i^*)^T (\nabla f(\bm{x}^{k}_i)-A^T_i\bm{\lambda}^* )+ \frac{M_i}{2}  \norm{\bm{x}^{k}_i-\bm{x}^{k+1}_i}^2. $$
	Recalling (\ref{first-order-updating}) that $\nabla f_i(\bm{x}^{k}_i) = A^T_i(\bm{\lambda}^k- \bm{\rho}^{k+1}_i(A_i\bm{x}^{k+1}_i+\sum_{j \not = i} A_j \bm{x}^{k}_j-\bm{c})) + \bm{\Gamma}^{k+1}_i(\bm{x}^{k}_i-\bm{x}^{k+1}_i)$, we have the following,
	\begin{equation}
	\label{expansion}
	\frac{m_i}{2} \norm{\bm{x}^{k+1}_i-\bm{x}_i^*}^2 \leq (\bm{x}^{k+1}_i-\bm{x}_i^*)^T ( A^T_i(\bm{\lambda}^k-\bm{\lambda}^*)- A^T_i\bm{\rho}^{k+1}_i \sum_{j=1}^{N} A_j (\bm{x}^{k}_j-\bm{x}^*_j) + \bm{D} _i(\bm{x}^{k}_i-\bm{x}^{k+1}_i))+ \frac{M_i}{2} \norm{\bm{x}^{k}_i-\bm{x}^{k+1}_i}^2. 
	\end{equation}
	Due to the approximation, we have a different bound as
	\begin{equation}
	\begin{aligned}
	\norm{\bm{\lambda}^{k+1}-\bm{\lambda}^*}^2 \leq& \frac{1}{\sigma^2_{i,\min} }\bigg\| \nabla f_i(\bm{x}^{k}_i) - \nabla f_i(\bm{x}^{*}_i) - A_i^T(\bm{\lambda}^{k}-\bm{\lambda}^{k+1}) - \bm{D}_i(\bm{x}^{k}_i-\bm{x}^{k+1}_i) + A_i^T\bm{\rho}^{k+1}_i\sum_{j=1}^{N}A_j( \bm{x}^{k}_j-\bm{x}^{*}_j)\bigg\|^2 \\
	\overset{(a)}{\leq}&  \frac{5}{\sigma^2_{i,\min} } ( 2M_i\norm{\bm{x}^{k+1}_i-\bm{x}^*_i }^2  + \sigma^2_{i,\max} \norm{\bm{\lambda}^{k}-\bm{\lambda}^{k+1}}^2 + (D^2_{i}+2M_i) \norm{ 
		\bm{x}^{k }_{i}-\bm{x}^{k+1}_i}^2 +  \\
	& \rho_{i,\max}^{2(k+1)}\sigma^2_{i,\max}\bigg\| \sum_{j=1}^{N}A_j( \bm{x}^{k}_j-\bm{x}^{k+1}_j)\bigg\|^2+\rho_{i,\max}^{2(k+1)}\sigma^2_{i,\max}\bigg\| \sum_{j=1}^{N}A_j( \bm{x}^{k+1}_j-\bm{x}^{*}_j)\bigg\|^2),
	\end{aligned}
	\end{equation}
	where $(a)$ is from the fact that $\norm{\nabla f_i(\bm{x}^{k}_i) - \nabla f_i(\bm{x}^{*}_i)}^2\leq M_i\norm{x_i^k-x_i^{k+1}+x_i^{k+1}-x_i^*}^2\leq 2M_i\norm{x^k_i-x^{k+1}_i}^2 + 2M_i\norm{x^{k+1}_i-x_i^*}^2$. The rest of the proof is similar to that of Theorem 3.1, and the $p$ can be selected as
	\begin{equation}
	\min \bigg\{ \frac{ m_i- \alpha  N\rho_{i,\max}^{2(k+1)}\sigma^2_{i,\max} - \alpha \breve{\rho}_{i,\max}^{2(k+1)}\sigma^2_{i,\max} }{\frac{10M_i}{\zeta N\sigma^2_{i,\min}}+D_i },\frac{D_i-\frac{N\sigma^2_{i,\max}}{\alpha }-M_i}{\frac{5 ( D^2_i+2M_i)}{\zeta N\sigma^2_{i,\min}} + \frac{5 \sigma^2_{i,\max}}{\zeta}\sum_{ j=1}^{N}\frac{\rho_{j,\max}^{2(k+1)}\sigma^2_{j,\max}}{\sigma^2_{j,\min}}} , \frac{\frac{2\rho^0}{\zeta }-\frac{N}{\alpha\zeta  }-1}{ 5  \sum_{i=1}^{N}(\frac{\rho_{i,\max}^{2(k+1)}}{\zeta^2}+\frac{1}{N}) \frac{\sigma^2_{i,\max}}{\sigma^2_{i,\min}}  } \bigg\},
	\end{equation}
	with parameters:
	\begin{equation}
	\label{rho-fa-selection}
	\left\{ 
	\begin{aligned}
	&\alpha < \frac{ m_i}{N\rho_{i,\max}^{2(k+1)}\sigma^2_{i,\max}+ \breve{\rho}_{i,\max}^{2(k+1)}\sigma^2_{i,\max}}, \\
	&D_i >\max\{ \rho_{i,\max}^{k+1}\sigma^2_{i,\max},\frac{N\sigma^2_{i,\max}}{\alpha}+M_i\}, \\
	&\rho^0 >\frac{N}{2\alpha},\\
	&\zeta < 2\rho^0 - \frac{N}{\alpha}.
	\end{aligned}
	\right.
	\end{equation}

	

	\section{Proof of Theorem \ref{noise-thm}}
	From the updating procedure with noise,
	\begin{equation}
	\bm{x} ^{k+1}_i = D_i^{-1} (A^T_i\bm{\rho}^{k+1}_i(\bm{c}-\sum_{j \not = i} A_j \bm{x}^{k}_j) + A_i^T \bm{\lambda }^{k} +\bm{\Gamma}^{k+1}_i \bm{x} ^{k}_i - \nabla f_i( \bm{x} ^k_i))+\bm{\Delta}^{k+1}_i.
	\end{equation}
	We then derive the expression of $\nabla f( \bm{x }^{k}_i)$ as follows,
	\begin{equation}
	\begin{aligned}
	\nabla f( \bm{x} ^{k}_i) = A^T_i \bm{\lambda }^k-A^T_i\bm{\rho}_i^{k+1}  \sum_{j=1}^{N} A_j(  \bm{x} ^{k}_j - \bm{x}_j^* ) + \bm{D}_i( \bm{x} ^{k }_i -  \bm{x} ^{k+1}_i) + \bm{D}_i\bm{\Delta}^{k+1}_i.
	\end{aligned}
	\end{equation}
	It is noted that the only difference, when compared to (\ref{gradient-1}), arises from the additional term $\bm{\Delta}^{k+1}_i$. Due to the strong convexity assumed, we conduct a similar reasoning as (\ref{expansion}) and have the following inequality: 
	\begin{equation}
	\begin{aligned}
	\label{noisy-basis}
	&\frac{m_i}{2} \norm{ \bm{x} ^{k+1}_i-\bm{x}_i^*}^2 \\
	& \leq ( \bm{x} ^{k+1}_i-\bm{x}_i^*)^T( A^T_i( \bm{\lambda} ^k-\bm{\lambda}^*)- A^T_i\bm{\rho}_i^{k+1}  \sum_{j=1}^{N} A_j( \bm{x} ^{k}_j - \bm{x}_j^* ) + \bm{D}_i( \bm{x} ^{k }_i -  \bm{x }^{k+1}_i) + \bm{D}_i\bm{\Delta}^{k+1}_i) + \frac{M_i}{2} \norm{ \bm{x} ^{k }_i- \bm{x} ^{k+1}_i}^2. 
	\end{aligned}
	\end{equation}
	By summing up over $i$ from 1 to $N$ on both sides of (\ref{noisy-basis}), we have 
	\begin{equation}
	\begin{aligned}
	\label{noisy-basis-1}
	&\sum_{i=1}^N \frac{m_i}{2} \norm{ \bm{x} ^{k+1}_i-\bm{x}_i^*}^2 \\
	&\leq  \sum_{i=1}^{N}( ( \bm{x} ^{k+1}_i-\bm{x}_i^*)^T( A^T_i( \bm{\lambda} ^k-\bm{\lambda}^*)- A^T_i\bm{\rho}_i^{k+1}  \sum_{j=1}^{N} A_j( \bm{x} ^{k}_j - \bm{x}_j^* ) + \bm{D}_i( \bm{x} ^{k }_i -  \bm{x }^{k+1}_i) + \bm{D}_i\bm{\Delta}^{k+1}_i) + \frac{M_i}{2} \norm{ \bm{x} ^{k }_i- \bm{x} ^{k+1}_i}^2). 
	\end{aligned}
	\end{equation}
	
	By moving the left hand side to the right hand side, and taking the term $\bm{D}_i\bm{\Delta}^{k+1}_i$ out of the summation, we have 
	\begin{equation}
	\label{noisy-basis-2}
	\begin{aligned}
	&  \underbrace{\sum_{i=1}^{N} (( \bm{x} ^{k+1}_i-\bm{x}_i^*)^T (A^T_i( \bm{\lambda }^k-\bm{\lambda}^*)- A^T_i\bm{\rho}_i^{k+1}  \sum_{j=1}^{N} A_j( \bm{x} ^{k}_j - \bm{x}_j^* ) + \bm{D}_i( \bm{x} ^{k }_i - \bm{x} ^{k+1}_i))}_{(1)}+ \\
	& \qquad\qquad\qquad\qquad\underbrace{ \frac{M_i}{2}  \norm{ \bm{x} ^{k }_i- \bm{x} ^{k+1}_i}^2-\frac{m_i}{2} \norm{ \bm{x} ^{k+1}_i-\bm{x}_i^*}^2)}_{(1)} +  \underbrace{\sum_{i=1}^{N}( \bm{x }^{k+1}_i-\bm{x}_i^*)^T\bm{D}_i\bm{\Delta}^{k+1}_i}_{(2)} \geq 0.
	\end{aligned}
	\end{equation}
	Therefore, the proof of Theorem \ref{first-order-converge} shown in Appendix J is an analysis on term (1). From Theorem \ref{first-order-converge}, there exists $p>0$ for parameters within the admissible range defined in (\ref{admissable-range}), $\norm{\bm{u}^{k} - \bm{u}^{*} }^2_G \geq (1+p) \norm{\bm{u}^{k+1} - \bm{u}^{*} }^2_G$. Now combining both terms (1) and (2) to show the convergence rate, it still holds with almost the same reasoning except one difference. Due to the noise, the upper bound of $\norm{\bm{\lambda}^{k+1}-\bm{\lambda}^{*}}^2$, given before as (\ref{upper-linear}), becomes 
	\begin{equation}
	\begin{aligned}
	&\norm{ \bm{\lambda} ^{k+1}- \bm{\lambda} ^{*}}^2 \\
	& \leq\frac{1}{\sigma^2_{i,\min} } \big\|\nabla f_i( \bm{x} ^{k}_i) - \nabla f_i(\bm{x}^{*}_i) - A_i^T( \bm{\lambda }^{k}- \bm{\lambda }^{k+1}) - \bm{D}_i( \bm{x} ^{k}_i- \bm{x} ^{k+1}_i) + A_i^T\bm{\rho}^{k+1}_i\sum_{j=1}^{N}A_j(  \bm{x} ^{k}_j-\bm{x}^{*}_j)+ \bm{D}_i\bm{\Delta}^{k+1}_i\big\|^2 \\
	& \leq \frac{6}{\sigma^2_{i,\min} } ( 2M_i\norm{ \bm{x} ^{k+1}_{i}-\bm{x}^{*}_{i} }^2 + \norm{A_{i}( \bm{\lambda} ^{k}- \bm{\lambda} ^{k+1})}^2 +(D^2_{i} +2M_i)\norm{  \bm{x} ^{k}_{i}- \bm{x} ^{k+1}_{i} }^2  +  \\
	&\rho_{i,\max}^{2(k+1)}\sigma^2_{i,\max}\big\| \sum_{j=1}^{N}A_j(  \bm{x }^{k}_j- \bm{x} ^{k+1}_j)\big\|^2+\rho_{i,\max}^{2(k+1)}\sigma^2_{i,\max}\big\|\sum_{j=1}^{N}A_j(  \bm{x} ^{k+1}_j-\bm{x}^{*}_j)\big\|^2+D_i^2\norm{\bm{\Delta}^{k+1}_i}^2).
	\end{aligned}
	\end{equation}
	The changes in the constants here slightly change the range of $p$ selection but do not affect the existence of $p$ such that 
	\begin{equation}
	\label{admm_up}
	\begin{aligned}
	\norm{\bm{u}^{k}-\bm{u}^{*}}^2_G   \geq& (1+p) \norm{\bm{u}^{k+1}-\bm{u}^{*}}^2_G - 2\sum_{i=1}^N D_i(\bm{x}^{k+1}_i-\bm{x}^*_i)^T \bm{\Delta}^{k+1}_i - \frac{6\delta}{\zeta N}  \sum_{i=1}^N\frac{D_i^2}{\sigma^2_{i,\min}} \norm{\bm{\Delta}^{k+1}_i}^2 \\
	\geq &(1+(1-\hat{\epsilon} )p) \norm{\bm{u}^{k+1}-\bm{u}^{*}}^2_G - \sum_{i=1}^N( \frac{6\delta D_i^2}{\zeta N\sigma^2_{i,\min}} + \frac{D_i}{\hat{\epsilon}p}) \norm{\bm{\Delta}^{k+1}_i}^2,
	\end{aligned}
	\end{equation}
	where $\hat{\epsilon}\in(0,1)$. Let $a = \frac{1}{1+(1-\hat{\epsilon})p}$, then \footnote{
		As a short comment, when $\lim_{K \to \infty} \norm{\bm{\Delta}^{K}_i}^2 \to 0$, there exists a constant $C$ that $\sum_{i=1}^{N}(\frac{6\delta D_i^2}{\zeta N\sigma^2_{i,\min}} +  \frac{D_i}{\hat{\epsilon}p})\norm{\bm{\Delta}^{K}_i}^2 \leq C \max_{i}\norm{\bm{\Delta}^{K}_i}^2$. Therefore,
		\begin{equation}
		R^{K+1} \leq  C\sum_{k=1}^{K+1} \max_i\norm{\bm{\Delta}^{k}_i}^2a^{K+2-k}. 
		\end{equation}
		For any arbitrarily small constant $z>0$, there exists $k_0$, such that for any $K>2k_0$, 
		$$ C\sum_{k=1}^{k_0} \max_i \norm{\bm{\Delta}^{k}_i}^2a^{K+1-k} \leq Ca^{k_0}\sum_{k=1}^{k_0} \max_i \norm{\bm{\Delta}^{K}_i}^2a^{k_0+1-k} <\frac{{z}}{2}.$$ On the other hand, $\max_i \norm{\bm{\Delta}^{k}_i}^2 \leq \frac{ z (1-c)}{2Cc}$ for any $k>k_0$. Therefore,
		\begin{equation}
		\begin{aligned}
		R^{K}\leq C\sum_{k=1}^{k_0} \max_i \norm{\bm{\Delta}^{k}_i}^2a^{K+1-k} + C\sum_{k=k_0+1}^{K} \max_i \norm{\bm{\Delta}^{k}_i}^2a^{K+1-k} \leq \frac{z}{2}+ C\max_i \big\|\bm{\Delta}^{k_0+1}_i\big\|^2 \sum_{k=k_0+1}^{K} a^{K-k} \leq  z.
		\end{aligned}
		\end{equation}
	}
	\begin{equation}
	\begin{aligned}
	\norm{\bm{u}^{K+1}-\bm{u}^{*}}^2_G \leq & a\norm{\bm{u}^{K}+\bm{u}^{*}}^2_G +  \sum_{i=1}^N( \frac{6p D_i^2}{\zeta N\sigma^2_{i,\min}} +  \frac{D_i}{\hat{\epsilon}p}) a\norm{\bm{\Delta}^{K+1}_i}^2\\
	\leq & \cdots\\
	\leq & a^{K+1}\norm{\bm{u}^{0}-\bm{u}^{*}}^2_G + \sum_{i=1}^N( \frac{6p D_i^2}{\zeta N\sigma^2_{i,\min}} +  \frac{D_i}{\hat{\epsilon}p}) \sum_{k=1}^{K+1}a^k\norm{\bm{\Delta}^{K+2-k}_i}^2\\ 
	= & a^{K+1}\norm{\bm{u}^{0}-\bm{u}^{*}}^2_G + R^{K+1}. 
	\end{aligned}
	\end{equation}

	At last, we analyze the utility-privacy tradeoff. With respect to the $\mathscr{B}_{\infty}$ sensitivity, based on the Laplace mechanism \cite{laplace}, let each coordinate of $\bm{\Delta}^k_i$ for any $k$, i.i.d. follow Lap$(0, \frac{\epsilon D_i }{\mathscr{B}_{\infty}Kd})$ for the composition across $K$ iterations and $d$ dimensions in pure $\epsilon$-LDP. Substituting the above form into (\ref{admm_up}), we have 
	$$ \norm{\bm{u}^{K}-\bm{u}^{*}}^2_G = O(a^K \norm{\bm{u}^{0}-\bm{u}^{*}}^2_G + 2Nd(\frac{\mathscr{B}_{\infty}Kd}{\epsilon})^2), $$
	which is $\tilde{O}(\frac{N\mathscr{B}^2_{\infty}d^3}{\epsilon^2})$ due to the exponential decaying of the first term. Here we omit all other constants to avoid the tedious expression on $\frac{1}{1-a}$. Similarly, under $(\epsilon,\delta)$-LDP, with the strong composition \cite{relaxedDP}, we only require that each coordinate of $\bm{\Delta}^k_i$ for any $k$, i.i.d. follow Lap$(0, O(\frac{\epsilon D_i}{\mathscr{B}_{\infty}\sqrt{Kd}}))$ and thus the utility loss is $\tilde{O}(\frac{N\mathscr{B}^2_{\infty}d^2}{\epsilon^2})$. \footnote{Note that in relaxed LDP, $\delta$ is assumed as a constant and $\epsilon$ is sufficiently small. Thus we drop the $\log(\frac{1}{\delta})$ term.}

	\section{Proof of Theorem 3.4}
	For each $i \in [1:N]$, 
	\begin{equation}
	\label{gd-convergence-1}
	\begin{aligned}
	&\norm{\bm{x}^{k+1}_i - \bm{x}^*}^2 = \norm{\frac{2\bm{w}^{k+1}_{i}}{N} \bm{x}^k_{i_k}+ \frac{2(\bm{I}_{d}-\bm{w}^{k+1}_{i})}{N} \bm{x}^k_{\hat{i}_k} + \frac{\sum_{j \not = i_k, \hat{i}_k} \bm{x}^k_j }{N} -\eta_{k+1}\nabla f_i(\bar{\bm{x}}^k )+\bm{\Delta}^{k+1}_i- \bm{x}^*}^2 \\
	& = \norm{\frac{2\bm{w}^{k+1}_{i}}{N} (\bm{x}^k_{i_k}-\bm{x}^*)+ \frac{2(\bm{I}_{d}-\bm{w}^{k+1}_{i})}{N} (\bm{x}^k_{\hat{i}_k}-\bm{x}^*) + \frac{\sum_{j \not = i_k, \hat{i}_k} (\bm{x}^k_j -\bm{x}^*)}{N}}^2 +  \eta_{k+1} ^2\norm {\nabla f_i(\bar{\bm{x}}^k)+   \eta_{k+1}^{-1} \bm{\Delta}^{k+1}_i  }^2\\
	& - 2\eta_{k+1} \langle \frac{2\bm{w}^{k+1}_{i}}{N} \bm{x}^k_{i_k}+ \frac{2(\bm{I}_{d}-\bm{w}^{k+1}_{i})}{N} \bm{x}^k_{\hat{i}_k} + \frac{\sum_{j \not = i_k, \hat{i}_k} \bm{x}^k_j }{N}-\bm{x}^* , \nabla f_i(\bar{\bm{x}}^k) +   \eta_{k+1}^{-1}\bm{\Delta}^{k+1}_i  \rangle,
	\end{aligned}
	\end{equation}
	where $\bar{\bm{x}}^k = \frac{1}{N}\sum_{i=1}^N\bm{x}_i^k$.
	In the following, we will use the following inequality that, if for $i \in [1:N]$, $\omega_i > 0$ and $\sum_{i=1}^{N} \omega_i=1$, then for arbitrary $N$ real numbers $y_{[1:N]}$, the following holds,
	\begin{equation}
	\label{fact_1}
	(\sum_{i=1}^N \omega_i y_i)^2 \leq \sum_{i=1}^{N} \omega_i y^2_i.
	\end{equation}
	It is noted that in (\ref{gd-convergence-1}), the sum of weights on $(\bm{x}^k_i-\bm{x}^*)$ is always the identity. With (\ref{fact_1}), by taking expectation on both sides of (\ref{gd-convergence-1}), we have 
	\begin{equation}
	\label{gd-convergence-2}
	\begin{aligned}
	\mathbb{E} [\norm{\bm{x}^{k+1}_i - \bm{x}^*}^2] \leq \sum_{i=1}^{N} \frac{\norm{\bm{x}^k_i-\bm{x}^*}^2}{N} + \eta_{k+1}^2(G^2 +  \mathbb{E}[(\eta_{k+1}^{-1}\bm{\Delta}^{k+1}_i )^2] ) - 2\eta_{k+1} \langle \bar{\bm{x}}^k -\bm{x}^*  , \nabla f_i(\bar{\bm{x}}^k ) \rangle.
	\end{aligned}
	\end{equation}
	By summing up both sides of (\ref{gd-convergence-2}) from $i=1,2,...,N$, we can bound $\sum_{i=1}^{N} f_i(\bar{\bm{x}}^k ) - f_i(x^*)$ as follows,
	\begin{equation}
	\label{gd-convergence-3}
	\begin{aligned}
	\sum_{i=1}^{N} f_i(\bar{\bm{x}}^k ) - f_i(x^*) & \leq \langle \bar{\bm{x}}^k  - \bm{x}^* , \sum_{i=1}^N \nabla f_i(\bar{\bm{x}}^k ) \rangle \\
	& \leq \eta_{k+1}^{-1} (\sum_{i=1}^{N} \norm{\bm{x}^k_i-\bm{x}^{*}}^2 -  \norm{\bm{x}^{k+1}_i-\bm{x}^{*}}^2) + \eta_{k+1}\sum_{i=1}^N (G^2+\mathbb{E}[(\eta_{k+1}^{-1}\bm{\Delta}^{k+1}_i )^2]).
	\end{aligned}
	\end{equation}
	The above is due to the fact that $\sum_{i=1}^N \nabla f_i(\bm{x}^*)= \bm{0}$ and we drop the constant $\frac{1}{2}$ for brevity. Before we can derive a global convergence analysis, we need to give an upper bound on $\norm{\bm{x}^{k+1}_i-\bm{x}^{*}}$ with the initial divergence $\norm{\bm{x}^0_i-\bm{x}^{*}}$ and the noise $\mathbb{E}[(\bm{\Delta}^{k+1}_i)^2]$. It is noted that, with rearrangement on (\ref{gd-convergence-3}) and the fact $\sum_{i=1}^{N} f_i(\bm{x}^k_i) - f_i(\bm{x}^*) \geq 0$, 
	\begin{equation}
	\begin{aligned}
	\label{domain-bound}
	\mathbb{E}[\sum_{i=1}^{N} \norm{\bm{x}^{k+1}_i-\bm{x}^{*}}^2] \leq \mathbb{E}[ \sum_{i=1}^{N} \norm{\bm{x}^{k}_i-\bm{x}^{*}}^2 + \eta_{k+1}^2 (G^2+( \eta_{k+1}^{-1}\bm{\Delta}^{k+1}_i)^2)].
	\end{aligned}
	\end{equation}
	When we select  $\eta_{k} = \frac{1}{c\sqrt{k}}$, $\sum_{k=1}^{K}\eta_{k}^2 = \sum_{k=1}^{N} \frac{1}{c^2k} \leq \frac{\log K +1}{c^2}$ since $k\leq K$. The above render an upper bound on $\norm{\bm{x}^{k}_i-\bm{x}^{*}}^2$ that 
	$$ \mathbb{E} [\sum_{i=1}^N \norm{\bm{x}^{k}_i-\bm{x}^{*}}^2] \leq \mathbb{E}[\sum_{i=1}^N \norm{\bm{x}^{0}_i-\bm{x}^{*}}^2 + \frac{\log K +1}{c^2}N(G^2 + V^2)], $$ with the assumption $\mathbb{E}[(\bm{\Delta}^{k}_i/\eta_{k})^2] \leq V^2$. Thus, (\ref{gd-convergence-3}) can be further formulated as
	\begin{equation}
	\begin{aligned}
	\label{gd-convergence-4}
	\sum_{k=0}^{K-1}\sum_{i=1}^{N} \frac{\mathbb{E} [f_i(\bar{\bm{x}}^k)] - f_i(x^*)}{K} \leq& \frac{\sum_{k=1}^{K-1}( \eta_{k+1}^{-1} -\eta_k^{-1})\sum_{i=1}^{N}\norm{ \bm{x}_i^k-\bm{x}^* }^2+ \eta_1^{-1} \sum_{i=1}^N\norm{ \bm{x}_i^0-\bm{x}^*}^2+ \sum_{k=0}^{K-1}\eta_{k+1}N(G^2+V^2)}{K}\\
	= & O(\frac{c\sqrt{K} \sum_{i=1}^{N}\norm{\bm{x}^0_i-\bm{x}^*}^2 +N c^{-1}(\log{K}+2)\sqrt{K+1}(G^2 +V^2)}{K}),
	\end{aligned}
	\end{equation}
	and $\mathbb{E}[\sum_{i=1}^N f_i(\sum_{k=0}^{K-1}\sum_{i=1}^{N}\bm{x}^k_i/NK) -f_i(\bm{x}^*)] \leq  \sum_{k=0}^{K-1}\sum_{i=1}^{N} \frac{\mathbb{E} [f_i(\bar{\bm{x}}^k)] - f_i(x^*)}{K}$. Here we use the trick of SGD proof that selecting such a sequence of decreasing step size. To finally disclose the utility-privacy tradeoff, we specify the parameter of noise. In pure $\epsilon$-DP setting, since the sensitivity is bounded by $\mathscr{B}_{\infty}$ in $l_{\infty}$, on each dimension we may add a noise following Lap$(0, \frac{\epsilon}{dK\eta_k\mathscr{B}_{\infty}})$ to produce a total $\epsilon$ loss from $d$ dimensions and $K$ iterations. Under the relaxed $(\epsilon, \delta)$-DP setting, with the strong composition theorem \cite{relaxedDP}, the variance $\mathbb{E}[(\bm{\Delta}^{k}_i/ \eta_k)^2]$ can be reduced to $\tilde{O}(d(\frac{\sqrt{dK}\mathscr{B}_{\infty}}{\epsilon})^2)$. Substituting those into (\ref{gd-convergence-4}), we complete the proof of the Theorem that 
	\begin{equation}
	\left\{ 
	\begin{aligned}
	& \text{pure}~ \epsilon-LDP:  \frac{c \sum_{i=1}^{N}\norm{\bm{x}^{0}_i-\bm{x}^{*}}^2}{\sqrt{K}} +\frac{N(\log K +2)(G^2 +\frac{d^3K^2\mathscr{B}_{\infty}^2}{\epsilon^2} )}{c\sqrt{K}}   \\
	&  ~~~~~~~~~~~~~~~~~~~~= \tilde{O} \bigg(\frac{ \sqrt{\sum_{i=1}^{N}\norm{\bm{x}^{0}_i-\bm{x}^{*}}^2} \sqrt{N}(G+\frac{d^{3/2}K\mathscr{B}_{\infty}}{\epsilon})}{\sqrt{K}}\bigg) \\
	& \text{relaxed}~ (\epsilon, \delta)-LDP: \tilde{O} \bigg(\frac{ \sqrt{\sum_{i=1}^{N}\norm{\bm{x}^{0}_i-\bm{x}^{*}}^2}\sqrt{N}(G+\frac{d\sqrt{K}\mathscr{B}_{\infty}}{\epsilon})}{\sqrt{K}}\bigg) = \tilde{O} \bigg(\frac{\sqrt{\sum_{i=1}^{N}\norm{\bm{x}^{0}_i-\bm{x}^{*}}^2}\sqrt{N}d\mathscr{B}_{\infty}}{\epsilon}\bigg),
	\end{aligned} 
	\right. 
	\end{equation}
	where $K$ is sufficiently large. 
	For the general communication graph, one may apply the "average process" technique used in \cite{GD2009}, which leads to a similar result.

	\section{Proof of Theorem 4.1}
	Here, we still assume the $N$ agents are connected in a complete graph for simplicity and redefine the updating subroutine for agent $i$ as 
	$$ \bm{x}^{k+1}_i = \bar{\bm{x}}^k - \eta_{k+1} \nabla f_i(\bar{\bm{x}}^k) + \bm{\Delta}^{k+1}_i. $$
	Let $\bar{\bm{\Delta}}^k = \frac{1}{N}\sum_{i=1}^{N}\bm{\Delta}_i^{k}$, $\bm{y}^{k+1}_i = \bar{\bm{x}}^k - \eta_{k+1} \nabla f_i(\bar{\bm{x}}^k)$ and $\bar{\bm{y}}^{k+1} = \bar{\bm{x}}^k-\frac{\eta_{k+1}}{N}  \sum_{i=1}^N \nabla f_i(\bar{\bm{x}}^k)$. Here, we fix $\eta_{k+1} = \frac{1}{M}$. It is noted that 
	\begin{equation}
	\label{subG}
	\begin{aligned}
	&f_i(\bar{\bm{x}}^{k+1}) - f_i(\bar{\bm{y}}^{k+1} ) =  \frac{1}{b_i} \sum_{j=1}^{b_i}\left[ \phi(\langle \bar{\bm{x}}^{k+1}, \bm{z}^i_j \rangle, l^i_j) - \phi(\langle \bar{\bm{y}}^{k+1}, \bm{z}^i_j \rangle, l^i_j)\right] \\
	& \leq \frac{L}{b_i} \sum_{j=1}^{b_i} \left|\big\langle \bar{\bm{x}}^{k+1} - \bar{\bm{y}}^{k+1}, \bm{z}^i_j\big \rangle \right|  = \frac{L}{b_i} \sum_{j=1}^{b_i} \left|\big\langle  \bar{\bm{\Delta}}^{k+1}, \bm{z}^i_j \big\rangle \right|. 
	\end{aligned}
	\end{equation}
	Under $\mathscr{B}_2$ sensitivity, we consider applying a Gaussian mechanism instead of using Laplace noise. With $\eta_k = \frac{1}{M}$, we assume $\bm{\Delta}^{k}_i$ for each $k$ and $i$ is i.i.d. following a Gaussian distribution $\mathscr{N}(\bm{0}, \tilde{O} (\frac{K\mathscr{B}^2_2}{M^2\epsilon^2})\times \bm{I}_d )$ for $(\epsilon,\delta)$-LDP \cite{FOCS2014}. Now we apply the following fact with respect to (sub)Gaussian variables \cite{subG}: if $\bm{b} = (b_1, b_2, ... ,b_d)$ is a vector of independent mean-zero $\sigma$-(sub)Gaussian, then $$\mathbb{E}[|\langle \bm{b}, \bm{c} \rangle|] \leq \sigma \norm{\bm{c}}$$ for $\bm{c} \in \mathscr{R}^d$. Thus, with the above assumptions, (\ref{subG}) can be further bounded as 
	\begin{equation}
	f_i(\bar{\bm{x}}^{k+1}) - f_i(\bar{\bm{y}}^{k+1} ) \leq \tilde{O}\bigg(\frac{C_2L\sqrt{K}\mathscr{B}_2}{M\epsilon} \bigg).
	\end{equation}
	
	On the other hand, with the Taylor theorem, 
	\begin{equation}
	\label{dsc}
	\begin{aligned}
	f_i(\bar{\bm{y}}^{k+1}) & \leq f_i(\bar{\bm{x}}^{k} ) - \frac{\eta_{k+1}}{N} \langle \sum_{j=1}^N \nabla f_j(\bar{\bm{x}}^{k}), \nabla f_i(\bar{\bm{x}}^{k}) \rangle + \frac{M}{2} \norm{ \frac{\eta_{k+1}}{N}  \sum_{j=1}^N \nabla f_j(\bar{\bm{x}}^{k}) }^2. \\
	\end{aligned}
	\end{equation}
	Let $F(\bm{x}) = \sum_{i=1}^{N} f_i(\bm{x})$. Since $f_i(\bm{x})$ is $m$-strongly convex, $F(\bm{x})$ is $Nm$-strongly convex. With the strong convexity, 
	$$ F(\bm{x}^*) -F(\bm{x}) \geq -\frac{1}{2Nm} \norm{\nabla F(\bm{x})}^2.  $$
	
	By summing up both sides of (\ref{dsc}) from $i=1$ to $N$ and subtracting $F(x^*)$, we have 
	\begin{equation}
	\label{dsc}
	\begin{aligned}
	F(\bar{\bm{y}}^{k+1})-F({x}^{*}) = \sum_{i=1}^N f_i(\bar{\bm{y}}^{k+1}) -f_i(x^*) &\leq F(\bar{\bm{x}}^{k} ) - F(\bm{x}^*) - \bigg(\frac{\eta_{k+1}}{N} -\frac{M \eta_{k+1} ^2}{2N}\bigg)  \norm{ \nabla F(\bar{\bm{x}}^{k})}^2 \\
	& = (1-\frac{m}{M})(F(\bar{\bm{x}}^{k} ) - F(\bm{x}^*)).
	\end{aligned}
	\end{equation}
	
	Putting things together, we have
	\begin{equation}
	\label{log}
	\begin{aligned}
	F(\bar{\bm{x}}^{k+1}) - F(\bm{x}^*) &\leq F(\bar{\bm{y}}^{k+1}) - F(\bm{x}^*) + \tilde{O}\bigg(\frac{NC_2L\sqrt{K}\mathscr{B}_2}{M\epsilon} \bigg) \\
	&\leq (1-\frac{m}{M})(F(\bar{\bm{x}}^{k} ) - F(\bm{x}^*)) + \tilde{O}\bigg(\frac{NC_2L\sqrt{K}\mathscr{B}_2}{M\epsilon} \bigg) \\
	&\leq (1-\frac{m}{M})^{k+1}(F(\bar{\bm{x}}^{0} ) - F(\bm{x}^*)) + \tilde{O}\bigg(\frac{M}{m}\frac{NC_2L\sqrt{K}\mathscr{B}_2}{M\epsilon} \bigg).
	\end{aligned}
	\end{equation}                                                  
	Therefore, $F(\bar{\bm{x}}^{K}) - F(\bm{x}^*) = \tilde{O}(\frac{\sqrt{M}NC_2L\mathscr{B}_2}{m^{3/2}\epsilon} )$, when $K = \tilde{O}(\frac{M}{m})$. \footnote{With such selection of $K$, the two terms summed in (\ref{log}) are in the same order with respect to $\tilde{O}$.} Here we use the fact that $e^{\log (1-\frac{m}{M})} \leq e^{-\frac{m}{M}}.$

	\section{Proof of Theorem 4.2}
	The proof of Theorem 4.2 is very similar to that of Theorem 4.1. When $\norm{\nabla f_i(\bm{x})}_1 \leq C_1$, 
	$$ \mathbb{E} \big[|f_i(\bm{x} + \bm{\Delta}) - f_i(\bm{x}) |\big] \leq \mathbb{E}\big[\sup_{\norm{\bm{c}}_1\leq C_1} |\langle \bm{c}, \bm{\Delta} \rangle | \big]$$
	Still from \cite{subG}, when $\bm{\Delta}$ is a vector of $d$ i.i.d $\sigma$-(sub)Gaussian variables, 
	\begin{equation}
	\label{l_1_control}
	\mathbb{E}\big[\sup_{\norm{\bm{c}}_1\leq C_1} |\langle \bm{c}, \bm{\Delta} \rangle | \big] \leq O\big(\sqrt{\log d}C_1\sigma\big).
	\end{equation}
	As for the noise, one can do a similar reasoning as that in Theorem 4.1.
	Under either $\mathscr{B}_1$ or $\mathscr{B}_2$ sensitivity, we may apply the Gaussian mechanism same as 4.1 and assume the noise $\bm{\Delta}^k_i$ following $\mathscr{N}(\bm{0},\tilde{O} (\frac{K\mathscr{B}^2_1}{M^2\epsilon^2})\times \bm{I}_d )$ or $\mathscr{N}(\bm{0},\tilde{O} (\frac{K\mathscr{B}^2_2}{M^2\epsilon^2})\times \bm{I}_d )$, respectively. With the error control as (\ref{l_1_control}), substituting the noise expression into (\ref{log}), the upper bound follows.

	In the following, we first introduce the random gradient descent in a private version. Without loss of generality, we only focus on the centralized case to optimize a function $f(\bm{x})$ and one can easily generalize the following results with the proof of Theorem 4.1 to the decentralized case. As assumed before, $f$ is $m$-strongly convex and $\hat{M}$-coordinate smooth. $\norm{\nabla f}_1 \leq C_1$.
	\begin{equation}
	\label{coordinate}
	\bm{x}^{k+1}[l_{k+1}] = \bm{x}^{k}[l_{k+1}] - \eta_{k+1} \nabla f(\bm{x}^{k})[l_{k+1}] + \Delta^{k+1}, 
	\end{equation}
	while $\bm{x}^{k+1}[l] = \bm{x}^{k}[l], l \not = l_{k+1}$ for $l=1,2,...,d$. Here, $l_{k+1}$ is i.i.d. uniformly selected from $[1:d]$. Intuitively, in(\ref{coordinate}), privacy loss is shrunk to a single dimension, while the step size also scales by $\frac{1}{d}$ compared to conventional GD. Here we assume that $\Delta^{k+1}$ is i.i.d. in a Gaussian distribution $\mathscr{N}(0,\sigma^2)$, where we will determine $\sigma$ soon. 
	
	In the following, we fix $\eta_{k+1} = 1/\hat{M}$. By Taylor series and the component Lipschitz assumption, at round $(k+1)$, 
	\begin{equation}
	\label{cd_pf}
	\begin{aligned}
	\mathbb{E}\big[f(\bm{x}^{k+1})\big]  & \overset{(a)}{\leq}\mathbb{E}\big[f(\bm{x}^{k}   - \eta_{k+1} \nabla f(\bm{x}^{k})[l_{k+1}] e_{l_{k+1}})\big] + \frac{C_1}{d} \mathbb{E}\big[|\Delta^{k+1}|\big] \\
	&\overset{(b)}{\leq}  f(\bm{x}^{k})  - \eta_{k+1}  \mathbb{E}\big[(\nabla f(\bm{x}^{k})[l_{k+1}])^2\big] + \frac{(\eta_{k+1})^2\hat{M}}{2} \mathbb{E}\big[(\nabla f(\bm{x}^{k})[l_{k+1}])^2\big]+\frac{C_1}{d} \sigma  \\
	& \overset{(c)}{\leq}   f(\bm{x}^{k}) -\frac{1}{2d\hat{M}}\norm{\nabla f(\bm{x}^k)}^2 + \frac{C_1}{d} \sigma.
	\end{aligned}
	\end{equation}
	Here, $(a)$ is due to the fact that subsampling on dimensions is uniform and in expectation, the coordinate Lipschitz of $f$ is bounded by $\frac{\mathcal{C}_1}{d}$. $(b)$ is the Taylor theorem on dimension $l_{k+1}$ and (c) is due to the fact that $\eta_{k+1} = 1/\hat{M}$ and $\mathbb{E}_{l_{k+1}} [(\nabla f(\bm{x}^{k})[l_{k+1}])^2] = \frac{\norm{\nabla f(\bm{x}^{k})}^2}{d}.$ 
	
	With the $m$-strongly convex assumption of $f(\bm{x})$, 
	$$  f(\bm{x}^*) - f(\bm{x}^k) \geq - \frac{1}{2m} \norm{\nabla f(\bm{x}^k)}^2 .$$
	Substituting the above into (\ref{cd_pf}), where $f(\bm{x}^*)$ is subtracted in both sides, we have 
	\begin{equation}
	\label{cd_pf_2}
	\begin{aligned}
	\mathbb{E}[f(\bm{x}^{k+1})]  - f(\bm{x}^*) & \leq (1-\frac{m}{d\hat{M}}) (\mathbb{E}[f(\bm{x}^{k})]  - f(\bm{x}^*)) + \frac{C_1}{d} \sigma\\
	& \leq (1-\frac{m}{d\hat{M}})^{k+1} (f(\bm{x}^{0})  - f(\bm{x}^*)) + \sum_{j=1}^{k+1} (1-\frac{m}{d\hat{M}})^{k+1-j} \frac{C_1}{d} \sigma \\
	& \leq e^{-\frac{m(k+1)}{d\hat{M}}} + \frac{d\hat{M}}{m} \cdot \frac{C_1}{d} \sigma. 
        \end{aligned}
	\end{equation}
	
	Now we turn to the privacy part. It is a bit different from before as we have to deal with the coordinate-wise sensitivity. We say $f(\cdot)$ is with $\mathscr{B}^c$ sensitivity if 
	$$\mathscr{B}^c \geq \sup_{\hat{f}, f \in \mathscr{F}} \sqrt{\sum_{l=1}^{d} \sup_{\bm{x} \in \mathcal{C}} \big(\nabla f(\bm{x})[l] - \nabla \hat{f}(\bm{x})[l]\big)^2}.$$ 
	$\mathscr{B}^c$ can be viewed as the $l_2$ norm of coordinate sensitivity. Recall the step size is $\frac{1}{\hat{M}}$, we consider fixing $\sigma^2 = \frac{(1+a)(\mathscr{B}^c)^2K}{\epsilon^2\hat{M}^2d}$ and we will determine $a$ in the following. First, from (\ref{log}), we select $K =\tilde{O}(d\frac{\hat{M}}{m})$. Across $K$ iterations, let $A_l$, $l \in [1:d]$, denote the number of times that $l$-th dimension selected as the descent direction. Due to the uniformness, $\mathbb{E}[A_l] =\mu= \frac{K}{d} = \tilde{O}(\frac{\hat{M}}{m})$. Now we apply Chernoff bound, 
	$$ \Pr(A_l\geq (1+a)\mu) \leq e^{-\frac{-a^2\mu}{2+a}}.$$
Therefore, for all $l \in [1:d]$, $A_l < (1+a)\mu$, the probability is upper bounded by 	$de^{-\frac{a^2\mu}{2+a}}$ with union bound. For any $f, \hat{f} \in \mathscr{F}$, let $\mathscr{B}^c_{l} = \sup_{\bm{x}} |\nabla f(\bm{x})[l] - \nabla \hat{f}(\bm{x})[l] |$. With the strong composition of privacy loss \cite{relaxedDP}, it is noted that
$$ \frac{(1+a)\epsilon^2(\mathscr{B}^c)^2K}{d\sum_{l=1}^{d} A_l (\mathscr{B}^c_{l})^2} >   \frac{\epsilon^2(\mathscr{B}^c)^2}{\sum_{l=1}^d (\mathscr{B}^c_{l})^2} \geq \epsilon^2 $$
when all $l \in [1:d]$, $A_l < (1+a)\mu$. For arbitrary constant failure probability, one may select $a =O( \frac{\log d}{\mu})$. To conclude, through RCD, one may obtain a utility bound in terms of $\hat{M}$ that $\tilde{O}(\frac{C_1\hat{M}\mathscr{B}^c}{m^2\epsilon}).$ with a similar reasoning.

\end{document}